\documentclass{amsart}
\usepackage{euscript,amssymb,xspace,mathrsfs}
\usepackage[author-year]{amsrefs}

\title[Brownian Sheet]{Images of the Brownian Sheet}%
   \thanks{The authors' research was
   supported by a grant from the National Science Foundation}
   \author[Khoshnevisan]{Davar Khoshnevisan}
   \address{Department\@ of Mathematics, The University of Utah,
      155 S.\@ 1400 E.\@ Salt Lake City, UT 84112--0090}
   \email{davar@math.utah.edu}
   \urladdr{http://www.math.utah.edu/\~{}davar}
\author[Xiao]{Yimin Xiao}
\address{Department of Statistics and Probability, A-413 Wells
    Hall, Michigan State University,
    East Lansing, MI 48824}
\email{xiao@stt.msu.edu}
\urladdr{http://www.stt.msu.edu/\~{}xiaoyimi}
\date{September 12, 2004}

\theoremstyle{plain}{
\newtheorem{theorem}{Theorem}[section]}
\theoremstyle{plain}{
\newtheorem{proposition}[theorem]{Proposition}}
\theoremstyle{plain}{
   \newtheorem{lemma}[theorem]{Lemma}}
\theoremstyle{plain}{
   \newtheorem{corollary}[theorem]{Corollary}}
\theoremstyle{definition}{
   }
\theoremstyle{definition}{
   }
\theoremstyle{remark}{
   \newtheorem{remark}[theorem]{Remark}}
\usepackage{amssymb,euscript}
\usepackage[author-year]{amsrefs}
\numberwithin{equation}{section}
\renewcommand{\O}{\varnothing}

\newcommand{\dimh}{\dim_{_\mathscr{H}}}

\newcommand{\dimp}{\dim_{_\mathscr{P}}}
\newcommand{\F}{\mathscr{F}}
\newcommand{\1}{\mathbf{1}}
\newcommand{\e}{\varepsilon}
\newcommand{\s}{\sigma}
\renewcommand{\P}{\mathrm{P}}
\newcommand{\E}{\mathrm{E}}

\newcommand{\C}{\mathrm{Cap}}
\newcommand{\R}{\mathbf{R}}
\newcommand{\Q}{\mathbf{Q}}
\renewcommand{\l}{\lambda}

\newcommand{\cle}{\preccurlyeq}
\newcommand{\lepi}{\preccurlyeq_\pi}
\newcommand{\supi}{\succcurlyeq_\pi}
\newcommand{\minpi}{\curlywedge_\pi}
\newcommand{\Fpi}{\F_\pi}
\begin{document}
\begin{abstract}
   An $N$-parameter Brownian sheet in $\R^d$
   maps a non-random compact set $F$
   in $\R^N_+$ to the random compact set $B(F)$
   in $\R^d$. We prove two results on the image-set $B(F)$:\\
   \indent{(1)} It has
   positive $d$-dimensional Lebesgue measure if
   and only if $F$ has positive $\frac d 2$-dimensional
   capacity. This generalizes greatly the earlier
   works of J.\@ Hawkes~\ycite{Hawkes},
   J.-P.\@ Kahane~\citelist{\ycite{Kahane85a}\ycite{Kahane85b}}, and one of
   the present authors~\ycite{Kh99}.\\
   \indent{(2)} If $\dimh  F > \frac d 2$, then
   with probability one, we can find a finite number of
   points $\zeta_1,\ldots,\zeta_m\in\R^d$ such that for
   any rotation matrix $\theta$ that leaves $F$
   in $\R^N_+$, one of the $\zeta_i$'s is interior to
   $B(\theta F)$. In particular,
   $B(F)$ has interior-points a.s.
   This verifies a conjecture of
   T.~S.~Mountford~\ycite{Mountford89}.\\
   \indent This paper contains two novel ideas:
   To prove {(1)}, we introduce
   and analyze a family of bridged sheets. Item
   {(2)} is proved by developing a notion
   of ``sectorial local-non-determinism (LND).''  Both ideas
   may be of independent interest.

   We showcase sectorial
   LND further by exhibiting some arithmetic properties of
   standard Brownian motion; this completes the work
   initiated by~\ocite{Mountford88}.
\end{abstract}
\maketitle\tableofcontents
\section{Introduction}
Let $B=\{B(t)\}_{t \in \R_+^N}$ denote the $(N,d)$-Brownian sheet.
That is, $B$ is the $N$-parameter Gaussian random field with values
in $\R^d$; its mean-function is zero, and its covariance function
is given by the following:
\begin{equation}
   \E\left[ B_i(s) B_j(t) \right] =
   \begin{cases}
      \prod_{k=1}^N \min(s_k,t_k),&\text{if $1\le i=j \le d$},\\
      0,&\text{otherwise}.
   \end{cases}
\end{equation}
We have written $B(t)$ in vector form as $(B_1(t),\ldots,B_d(t))$,
as is customary.

When $N=1$, $B$ is just Brownian motion in $\R^d$. In this case,
it is well known~\cite{Hawkes} that for any non-random compact set
$F\subseteq\R_+$,
\begin{equation}\label{eq:main}
   \P\left\{ \l_d(B(F))>0\right\}>0\ \quad\text{if and only if}
   \quad\ \C_{d/2}(F)>0.
\end{equation}
Here, $\l_d$ denotes the $d$-dimensional Lebesgue measure, and for
all $\alpha>0$, $\C_\alpha(F)$ denotes the $\alpha$-dimensional
Bessel--Riesz capacity of $F$ based on the $\alpha$-dimensional
energy form $I_\alpha$; i.e.,
\begin{equation}
   \C_\alpha(F) = \left[ \inf_{\mu\in\mathscr{P}(F)}
   I_\alpha(\mu)\right]^{-1},
   \text{ where }I_\alpha(\mu)
   =\iint \frac{\mu(ds)\, \mu(dt)}{\|s-t\|^\alpha},
\end{equation}
and $\mathscr{P}(F)$ denotes the collection of all probability
measures that are supported in $F$.

According to Taylor's theorem~\cite{Kh}*{Corollary 2.3.1, p.\@
525}, for all $F\subset\R^N_+$, $\C_\alpha(F)=0$ except possibly
when $\alpha<N$. Therefore, when $N=1$, (\ref{eq:main}) has
nontrivial content when, and only when, $d=1$.

In order to go beyond the one-dimensional case,~\ocite{Kahane85a}
proposed considering $N$-parameter processes (i.e., fractional
Brownian motion), and devised a Fourier-analytic argument which,
in the present setting, implies the following for the Brownian
sheet:
\begin{equation}
   \C_{d/2}(F)>0\ \Longrightarrow\
   \P\left\{ \l_d(B(F))>0\right\}>0\
   \Longrightarrow\
   \mathrm{H}_{d/2}(F)>0.
\end{equation}
Here, $\mathrm{H}_\alpha$ denotes the $\alpha$-dimensional
Hausdorff measure~\cite{Kahane85a}*{p.\@ 131, Remark 4}. There is
an obvious gap between the enveloping conditions of positive
capacity and measure. In the special case that $N=2$, this gap was
closed in~\ocite{Kh99}, but the problem for $N>2$ has remained
open. One of the intentions of this article is to complete the
existing picture by deriving the following:

\begin{theorem}\label{th:Kahane}
   For any choice of $N$ and $d$, and for all non-random compact sets
   $F$, the $d$-dimensional Lebesgue measure of $B(F)$ is positive
   with positive probability if and only if $\C_{d/2}(F)>0$.
\end{theorem}

We will prove also that the following is an equivalent formulation
of Theorem~\ref{th:Kahane}.

\begin{theorem}\label{th:level}
   For any choice of $N$ and $d$, for all non-random compact sets $F$,
   and  for every $a \in \R^d$,
   $\P\{B^{-1}(\{a\}) \cap F \neq\O \} > 0$ if and
   only if $\C_{d/2}(F)>0$.
\end{theorem}

On one hand, this and Taylor's theorem together show that if $N\le
\frac d2$, then $B^{-1}(\{x\})=\O$ almost surely for all
$x\in\R^d$. On the other hand, when $N>\frac d2$, the codimension
of $B^{-1}(\{x\})$ is almost surely $\frac d2$~\cite{Kh}*{\S4.7,
p.\@ 435}. In particular, we can conclude from Theorem 4.7.1
of~\ocite{Kh}*{p.\@ 436} that the Hausdorff--Besicovitch dimension
of $B^{-1}(\{x\})$ is a.s.\@ $N-\frac d2$. When $d=1$, this last
assertion is due to
Adler~\citelist{\ycite{Adler}\ycite{Adler'}\ycite{Adler0}}. The
general case $1\le d<2N$ was treated by~\ocite{Ehm}.

In fact, one can go a bit farther at little extra cost. Suppose
$f:\R_+\to\R_+\cup\{\infty\}$ is a non-increasing measurable
function that is finite everywhere except possibly at zero. We can
then define the \emph{$f$-capacity} of a
Borel set $F\subseteq\R^N_+$ as
\begin{equation}
   \C_f (F) = \left[ \inf_{\mu\in\mathscr{P}(F)}
   I_f (\mu)\right]^{-1},
   \text{ where }I_f(\mu)
   =\iint f(\|s-t\|)\, \mu(ds)\, \mu(dt).
\end{equation}
After combining our Theorem~\ref{th:level} with Theorem 15.2
of~\ocite{Peres}, we immediately obtain the following extension of
Theorem 5 of~\ocite{Hawkes}.

\begin{corollary}
   Let $f:\R_+\to\R_+\cup\{\infty\}$ be a non-increasing
   measurable function that is finite on $(0,\infty)$.
   Then for all $a\in\R^d$,
   \begin{equation}
      \P\left\{ \C_f\left( B^{-1}(\{a\}) \right)>0\right\}>0\
      \Longleftrightarrow\ \int_0^1
      \frac{f(t)\, dt}{t^{(d/2) -N+1}}<\infty.
   \end{equation}
\end{corollary}

Our proof of Theorem~\ref{th:Kahane} depends on: (i) Ideas from
the potential-theory of multiparameter processes that are nowadays
considered standard; and (ii) a novel analysis of a class of
embedded bridged sheets. We write Theorem~\ref{th:Kahane}
to not only document it in its definitive form, but to also
highlight some of the features of the said bridges. This
bridge-analysis is used in our forthcoming paper with Robert
Dalang and Eulalia Nualart to solve an old open problem on the
self-intersections of Brownian sheets.

Thanks to the Frostman theorem of potential theory and
Theorem~\ref{th:Kahane}, $\dimh  F>\frac d2$ implies that $B(F)$
can have positive Lebesgue measure, whereas $\dimh F<\frac d2$
implies that $B(F)$ almost-surely has zero Lebesgue measure. We
plan to prove that much more is true: ``\emph{If $\dimh  F > \frac
d2$, then $B(F)$ has interior-points almost surely.}''

This type of interior-point problem was first studied
by~\ocite{Kaufman} in the case of one-dimensional Brownian motion
($N=1$). In this case, Kaufman proved that if $\dimh  F >
\frac12$, then $B(F)$ has interior-points a.s.

Kahane~\citelist{\ycite{Kahane85a}\ycite{Kahane85b}}
and~\ocite{Pitt78} have extended Kaufman's result to symmetric
stable L\'evy processes and fractional Brownian motion,
respectively.

Mountford~\ycite{Mountford89} has considered such interior-point
problems for the Brownian sheet, and proved that if $\dimh  F >
\frac d2$, then for almost every rotation
$\theta F$ of $F$ that is in $\R^N_+$,
$B(\theta F)$ has interior-points
a.s.\footnote{``Almost every rotation''
holds with respect to the Haar measure on
rotation matrices.} Moreover, he has conjectured that $B(F)$ has
interior-points a.s.~\cite{Mountford89}*{p.\@ 184}. We verify
this conjecture by proving that the Brownian sheet has the
following striking property:

\begin{theorem}\label{th:interior}
   Let $B$ denote the $(N, d)$-Brownian sheet,
   and let $F \subset \R_+^N$
   be any non-random Borel set that satisfies
   $\dimh  F > \frac d2$.
   Then there a.s.\@ exist
   $\zeta_1,\ldots,\zeta_m\in\R^d$
   such that for every rotation matrix $\theta$
   that leaves $\theta F$ in $\R^N_+$
   we can find $1\le j\le m$ such that $\zeta_j$ is
   interior to $B(\theta F)$. In particular, let
   $\theta$ denote the identity
   to see that $B(F)$ has interior-points a.s.
\end{theorem}

Although fractional Brownian motion is
locally non-deterministic, the Brownian sheet is not. This
remark accounts for the differences between the methods
of~\ocite{Pitt78} and~\ocite{Mountford89}. As
part of our arguments, we prove that the Brownian sheet
satisfies a type of ``sectorial local-non-determinism''
(Proposition \ref{Prop:LND}); this property leads to a unification
of many of the methods developed for the fractional Brownian
motion and those for the Brownian sheet. We will show this,
anecdotally, by describing an improvement to older results
of~\ocite{Mountford88} on self-intersections of images of ordinary
Brownian motion.

The rest of this paper is organized as follows. Sections 2 and 3
reviews briefly the order structure of $\R^N$ and the commuting
property of the filtrations associated to the Brownian sheet.
Sections 4 and 5 describe the correlation structure of the
Brownian sheet, sectorial local non-determinism, and an
a class of bridged sheets. Theorems~\ref{th:Kahane} and
\ref{th:level} are proved in Sections 6 and 7, respectively. In
Section 8, we prove Theorem~\ref{th:interior}. We make further
remarks on the images of Brownian motion, and more general Gaussian
random fields, in Sections 9 and 10.

Unspecified positive and finite constants are denoted by $A$.
They are usually numbered by the equation in which they
appear.

\section{The Order Structure of $\R^N$}

We need to introduce a good deal of notation
in order to exploit the various Markov properties of
$B$ ``in various directions.''
This is the sole task of the present section.

\subsection{The Partial Orders}
There are $2^N$ natural partial orders on $\R^N$. There is a
convenient way to represent them all. Define,
\begin{equation}\label{eq:Pi}
   \Pi_N = \text{The power set of }\left\{ 1,\ldots,N\right\}.
\end{equation}
Then, each $\pi\in\Pi_N$ can be identified with the partial order
$\lepi$ on $\R^N$ as follows: For all $a,b\in\R^N$,
\begin{equation}
    a \lepi b\text{ iff for all } i=1,\ldots,N,\
    \begin{cases}
       a_i\le b_i, & \text{ if $i\in \pi$},\\
       a_i\ge b_i, & \text{ if $i\not\in \pi$}.
    \end{cases}
\end{equation}
We always write $\cle$ in place of the more
cumbersome $\cle_{\{1,\ldots,N\}}$.

An important feature of the totality $\{ \lepi\}_{\pi\in\Pi_N}$ of
these partial orders is that together they order $\R^N$. By this
we mean that for all $a,b\in \R^N$, there exists $\pi =
\pi(a,b)\in\Pi_N$ such that $a\lepi b$. [Simply, let
$\pi(a,b)=\{1\le i\le N:\ a_i\le b_i\}$.]

\subsection{The PO-Minimum}
Each partial order $\lepi$ naturally yields a $\pi$-minimum
operation $\minpi$ which we describe next.

For each point $b\in\R^N$, define $S_b^\pi$ to be its ``shadow in
the direction $\pi$''; i.e.,
\begin{equation}
   S_b^\pi = \left\{ a\in\R^N:\ a \lepi b \right\}.
\end{equation}
Then, given $a,b\in\R^N$ and a partial order $\pi\in\Pi_N$, we
define $a\minpi b$ to be the unique point whose shadow in the
direction $\pi$ is precisely $S_a^\pi\cap S_b^\pi$. Let us
emphasize the fact that
\begin{equation}
   c\lepi a ~,~ c\lepi b\ \Longrightarrow\
   c\lepi (a\minpi b).
\end{equation}
It is easy to prove that such a point always exists.

Each partial order $\pi\in\Pi_N$ on $\R^N$ induces $N$ linear
orders $\preccurlyeq_{\scriptscriptstyle (\pi,1)},
\cdots,\preccurlyeq_{\scriptscriptstyle (\pi,N)}$ on $\R$  via the
following:
\begin{equation}
    \preccurlyeq_{\scriptscriptstyle (\pi,\ell)} =
    \begin{cases}
       \le , & \text{ if $\ell \in \pi$},\\
       \ge , &\text{  if $\ell \not\in \pi$},
    \end{cases} \qquad {}^\forall \ell=1,\ldots, N.
\end{equation}
Of course, one obtains only two distinct partial orders this way:
$\le$ and $\ge$. However, in what is to come, the preceding
notation will seemlessly do most of the book-keeping for us.

%\section{Correlation Structure of Brownian Sheet}

\section{The Associated Filtrations}

Consider the $\s$-algebras
\begin{equation}
    \Fpi(t) =  \s\left\{ B(s) \right\}_{s\lepi t},\qquad{}^\forall
    t\in\R^N_+,\,  \pi\in\Pi_N.
\end{equation}
Informally speaking, knowing $\Fpi(t)$ amounts to knowing the
portion of the Brownian sheet $B$ that corresponds to the values
of $s$ in $\R^N_+$ that are less than $t$ in the partial order
$\pi$.

It is not difficult to see that for each partial order
$\pi\in\Pi_N$, the collection $\Fpi = \{ \Fpi(t)\}_{t\in\R^N_+}$
is a filtration indexed by $(\R^N_+,\lepi)$; i.e.,
\begin{equation}
    s\lepi t \ \Longrightarrow\
    \Fpi(s) \subseteq  \Fpi(t).
\end{equation}

For each partial order $\pi\in\Pi_N$, we also define $N$ {\em
one-parameter} families of $\s$-algebras, $\F^1_\pi,\ldots,
\F^N_\pi$, as follows:
\begin{equation}
    \F_\pi^\ell(r) = \sigma\left(
    B(s) ; ~ s_\ell \preccurlyeq_{\scriptscriptstyle (\pi,\ell)} r \right),
    \qquad\qquad  {}^\forall r\in\R_+.
\end{equation}
Note that $\F^\ell_\pi$ is a filtration of $\s$-algebras indexed
by $(\R_+,\preccurlyeq_{\scriptscriptstyle (\pi,\ell)})$.
Moreover, for all $t\in\R^N_+$, $\Fpi(t) = \cap_{\ell=1}^N
\F^\ell_\pi(t_\ell)$.

Following~\ocite{Kh}*{Chapter 1}, we say that $\Fpi$ is {\em
commuting}, if for all times $t\in\R^N_+$, the $\s$-algebras
$\F^1_\pi(t_1),\ldots, \F^N_\pi(t_N)$ are conditionally
independent given $\Fpi(t)$. This is a slightly more general
``F4-type'' property than the one of~\ocite{CW}.

\begin{proposition}\label{pr:commutation}
    For every $\pi\in\Pi_N$,
    the filtration $\Fpi$ is commuting in the partial order
    $\pi$; i.e., for all bounded
    random variables $Z$,
    \begin{equation}
        \E \left[\left. Z\, \right|\, \Fpi( s \minpi t) \right] =
        \E \left(
        \E \left[\left. Z\, \right|\, \Fpi(s) \right] \, \Big|\,
        \Fpi(t) \right),\qquad\textnormal{a.s.}
    \end{equation}
\end{proposition}
Thus, commuting filtrations refers to the commutation of the
conditional expectation operators.

\begin{proof}
   Define, for all $\pi\in\Pi_N$ and $t\in(0,\infty)^N$, define
   $\mathscr{I}(t)\in(0,\infty)^N$ coordinatewise as follows:
   \begin{equation}
      {\mathscr{I}}_j(t) =
      \begin{cases}
         t_j, &\text{ if $j\in \pi$},\\
         1/t_j, &\text{ if $j\not\in \pi$}.
      \end{cases}
   \end{equation}
   One can think of the map ${\mathscr{I}}$ as
   ``inversion off of $\pi$.''

   Now consider the following stochastic process,
   \begin{equation}
      W_\pi(t) = \frac{B \left( {\mathscr{I}}(t) \right)}{%
      \prod_{j\not\in \pi} {\mathscr{I}}_j (t) },
      \qquad{}^\forall t\in\R^N_+.
   \end{equation}
   This is a  Brownian sheet, as can be checked by
   computing covariances. Moreover,
   \begin{equation}
      \s\left( \left\{ W_\pi(s);\
      s\cle t \right\}\right) = \s\left( \left\{
      B\left( {\mathscr{I}}(s) \right);\ s\cle t\right\} \right)
      =\Fpi(t).
   \end{equation}
   Because the filtration, in the partial order $\cle$, of Brownian sheet
   is commuting~\cite{Kh}*{Theorem 2.4.1, p.\@ 237}, this shows
   that $\Fpi$ is also commuting.
\end{proof}

The preceding leads us to the following useful representation.

\begin{corollary}\label{co:commutation}
    For every $\pi\in\Pi_N$, $j=1,\ldots, N$,
    and $r \in \R_+$, define
    the conditional expectation operator,
    $\mathscr{E}^{j,r}_{_{\pi}} Y = \E[ Y\, |\,
    \F^j_{_{\pi}}(r) ]$ $(Y\in L^1(\P))$. Then,
    for all $t\in\R^N_+$ and for all $\P$-intergable
    random variables $Z$,
    \begin{equation}
        \E\left[ Z\, \left|\, \Fpi(t) \right. \right] =
        \mathscr{E}^{1,t_1'}_{_{\pi}}\cdots
        \mathscr{E}^{N,t_N'}_{_{\pi}} Z,
    \end{equation}
    where $(t_1',\ldots,t_N')$ denotes an arbitrary non-random permutation
of
    $t = (t_1,\ldots,t_N)$. Thus,
    \begin{equation}\label{eq:Cairoli}
       \E\left[ \sup_{t\in\Q^N_+} \left( \E\left[ Z\,
       \left|\, \Fpi(t) \right. \right] \right)^2 \right] \le
       4^N \E\left[ Z^2 \right].
    \end{equation}
\end{corollary}

\begin{proof}
   To prove the first display, we simply
   follow along the proof of Theorem 3.6.1 of~\ocite{Kh}*{p.\@ 38}, but
   everywhere replace
   $\cle$ and $\curlywedge$ by $\lepi$ and $\curlywedge_\pi$,
   respectively. For the second portion,
   we apply Doob's strong $(p,p)$-inequality
   for ordinary martingales $N$ times in succession. For
   example, see the proof of
   Cairoli's strong $(p,p)$-inequality~\cite{Kh}*{Theorem 2.3.1, p.\@ 19},
   but replace $\cle$ everywhere by $\lepi$.
\end{proof}

\section{Sectorial Local-Nondeterminism} \label{sec:LND}

In this and the next section we state and prove some results on the
correlation structure of the Brownian sheet $B =\{B(t)\}_{t \in
\R_+^N}$ in $\R^d$. In particular, we prove that $B$ is
sectorially locally non-deterministic, and that there is a natural
class of bridged sheets associated to $B$. These properties will
play an important role in this paper, as well as in studying the
self-intersections of the Brownian sheet.\\

\noindent\textbf{Assumption}\
Throughout Sections \ref{sec:LND} and \ref{Sec:Bridge}, we assume
that $d=1$.\\

The following lemma is well known; cf.\@ Lemmas 8.9.1 and 8.9.2 of
\ocite{Adler0}. For the sake of completeness, we describe a
simpler proof.
\begin{lemma}\label{lem:Q(0)}
   Choose and fix two numbers $0<a<b<\infty$.
   If $u,v\in[a,b]^N$, then
   \begin{equation}
      \frac{a^{N-1}}{\sqrt N} \left\| u-v\right\|\le
      \E\left[\left( B(u)-B(v)\right)^2 \right] \le N
      b^{N-1}\left\| u-v\right\|.
   \end{equation}
\end{lemma}

\begin{proof}
   Let $\s(s,t)=\prod_{j=1}^N s_j - \prod_{j=1}^N t_j$,
   and define $s\curlywedge t$ to be the vector whose $i$th coordinate
   is $s_i\wedge t_i$. Then clearly,
   \begin{equation}
      \E\left[\left( B(u)-B(v)\right)^2 \right]
      = \s\left( u,u \curlywedge v\right) +\s\left(
      v,u\curlywedge v\right).
   \end{equation}
   Clearly,
   \begin{equation}
      a^{N-1} \max_{1\le j\le N}(u_j-u_j\wedge v_j)\le
      \s\left(u,u\curlywedge v\right)\le N
      b^{N-1} \max_{1\le j\le N} (u_j-u_j\wedge v_j).
   \end{equation}
   A similar expression holds for $\s(v,u\curlywedge v)$, but everywhere
   replace $u_j-u_j\wedge v_j$ with $v_j-u_j\wedge v_j$. Add the two series
   of inequalities, and use the fact that
   $u_j+v_j-2(u_j\wedge v_j)=|u_j-v_j|$, to obtain
   \begin{equation}
      a^{N-1}\max_{1\le j\le N}|u_j-v_j| \le
      \E\left[\left( B(u)-B(v)\right)^2 \right] \le N
      b^{N-1}\max_{1\le j\le N}|u_j-v_j|.
   \end{equation}
   The lemma follows from this and the elementary fact
   that for all $N$-vectors ${x}$,
   $N^{-1/2}\| {x}\|\le
   \max_{1\le j\le N}|x_j|\le\| {x}\|.$
\end{proof}

The Brownian sheet is not locally non-deterministic (LND) with
respect to the incremental variance $\E\left[\left(
B(u)-B(v)\right)^2 \right]$. However, it satisfies the following
``sectorial'' type of local non-determinism; cf.\@
\ocite{Kh}*{Lemma 3.3.2, p.\@ 486} for a prefatory version.

\begin{proposition}[Sectorial LND]\label{Prop:LND}
   For all positive real number $a$,
   integers $n \ge 1$, and all
   $u,v,t^1, \ldots, t^n \in [a,\infty)^N$,
   \begin{align}
      &{\rm Var} \left( \left.  B(u)\, \right|\, B(t^1), \ldots ,
         B(t^n) \right) \ge \frac{a^{N-1}}{2}\
         \sum_{k=1}^N \min_{1\le j\le n}
         \left| u_k - t^j_k \right|,\label{Eq:LND}\\
      \begin{split}
         &{\rm Var} \left( \left.  B(u)-B(v)\, \right|\, B(t^1), \ldots ,
            B(t^n) \right)\\
         &\ \ge \frac{a^{N-1}}{2}\
            \sum_{k=1}^N \min\left(
            \min_{1\le j\le n}
            \left| u_k - t^j_k \right|+\min_{1\le j\le n}
            \left| v_k - t^j_k \right| ~,~ |u_k-v_k|\right).\label{Eq:LND2}
      \end{split}
   \end{align}
\end{proposition}

The proof is divided in two distinct steps. The first is the
analysis of the $N=1$ case; we present this portion next.

\begin{lemma}\label{lem:LND}
   Let $\{X(t)\}_{t\ge 0}$ denote standard Brownian motion
   on the line. Then for all times $s,t,s_1,\ldots,s_m\ge 0$,
   \begin{align}
      \mathrm{Var} \left( \left. X(s)\, \right|\,
         \mathscr{X} \right) & \ge \frac12
         \min_{1\le j\le m} \left| s-s_j \right|,\label{eq:BM1}\\
      \mathrm{Var} \left( \left. X(t)-X(s)\, \right|\,
          \mathscr{X} \right) & \ge \frac12
         \min\left( \min_{1\le j\le m} \left|s-s_j \right|
         + \min_{1\le j\le m} \left|t-s_j \right|
         ~,~ |t-s|\right), \label{eq:BM2}
   \end{align}
   where $\mathscr{X}$ denotes the $\s$-algebra generated
   by $(X(s_1),\ldots,X(s_m))$.
\end{lemma}

\begin{proof}
   Equation (\ref{eq:BM1}) follows from (\ref{eq:BM2}). Indeed,
   let $t=s_j$ in (\ref{eq:BM2}), and then optimize over
   all $j$ to obtain (\ref{eq:BM1}).

   Equation (\ref{eq:BM2}) is proved by analyzing two different
   cases. Throughout, we assume,  without any
   loss of generality, that $s<t$.\\
   \indent\emph{Case 1:}\ The first case is where some $s_j$ falls between
   $s$ and $t$. Recall that if $\F$ and $\F'$
   are linear subspaces (equivalently, $\s$-algebras)
   in the Gauss space $L^2(\P)$, then for every
   Gaussian variate $G\in L^2(\P)$,
   \begin{equation}\label{eq:projection}
      \F\subset\F'\ \Longrightarrow\
      \mathrm{Var}(G\,|\,\F)
      \ge\mathrm{Var}(G\,|\,\F').
   \end{equation}
   Moreover, both conditional variances are non-random.
   This elementary fact, used in conjunction with the Markov property,
   allows us to assume without any further loss in generality
   that $m=4$ and $s_1<s<s_2<s_3<t<s_4$.

   Now define
   $\xi_1=X(s_1)$, $\xi_2=X(s)-X(s_1)$, $\xi_3=X(s_2)-X(s)$,
   $\xi_4=X(s_3)-X(s_2)$, $\xi_5=X(t)-X(s_3)$, and $\xi_6=
   X(s_4)-X(t)$. These are independent Gaussian variables,
   and $\mathscr{X}$ is the linear subspace of $L^2(\P)$ that is spanned
   by $\xi_1$, $(\xi_1+\xi_2+\xi_3)$, $(\xi_1+\xi_2+\xi_3+\xi_4)$,
   and $(\xi_1+\xi_2+\xi_3+\xi_4+\xi_5+\xi_6)$.
   Therefore, by the independence of the $\xi_j$'s,
   \begin{equation}\begin{split}
      &\mathrm{Var}\left( \left.
         X(t)-X(s)\, \right|\, \mathscr{X}\right)\\
      & = \mathrm{Var}\left( \left. \xi_3+\xi_5\,\right|\,
         \xi_2+\xi_3, \xi_4, \xi_5+\xi_6\right)\\
      & = \inf_{\alpha,\beta\in\R}\E\left[\left(
         \xi_3+\xi_5 -\alpha(\xi_2+\xi_3) -
         \beta(\xi_5+\xi_6)\right)^2\right]\\
      & =\frac{(s_2-s)(s-s_1)}{s_2-s_1} +
         \frac{(s_4-t)(t-s_3)}{s_4-s_3},
   \end{split}\end{equation}
   whence (\ref{eq:BM2}) in the present case.\\
   \indent\emph{Case 2:}\ The remaining case is where no $s_j$ falls
   in $(s,t)$. In this case, the Markov property shows that
   we can assume, without loss of generality, that $m=2$
   and $s_1 < s < t< s_2$. A direct calculation reveals that
   in this case,
   \begin{equation}
      \mathrm{Var}\left( \left. X(t)-X(s)\,\right|\, \mathscr{X}\right)
      = \frac{(t-s)(s_2-s_1-t+s)}{s_2-s_1},
   \end{equation}
   from which (\ref{eq:BM2}) follows.
\end{proof}

\begin{proof}[Proof of Proposition~\ref{Prop:LND}]
   Let $\langle a \rangle = (a, \ldots, a)$
   designate the lower-left corner of $[a, \infty)^N$,  and for all $r\ge 0$
   and $1 \le k \le N$ define
   \begin{equation}
      X_k(r) = \frac{B(
      \overbrace{a, \ldots, a}\limits^{\text{$k-1$ terms}},
      a+r, a, \ldots, a) - B(\langle a
      \rangle)}{a^{(N-1)/2}}.
   \end{equation}
   The process $\{ X_k(r)\}_{r\ge 0}$ is a standard Brownian motion on the
   line.

   For all $t\in[a,\infty)^N$, we decompose the rectangle $[0, t]$ into the
   following disjoint union:
\begin{equation}\label{Eq:Split0}
[0, t] = [0, a]^N \cup \bigcup_{j=1}^N D(t_j) \cup \Delta(a, t),
\end{equation}
where $D(t_j) =\{ s \in [0, 1]^N: 0 \le s_i \le a \hbox{ if } i
\ne j,\, a < s_j \le t_j\}$ and $\Delta(a, t)$ can be written as a
union of $2^N - N -1$ sub-rectangles of $[0, t]$. Then we have the
following decomposition:  For all $t\in[a,\infty)^N$,
   \begin{equation}\label{Eq:decomp}
      B(t) = B(\langle a \rangle) + a^{(N-1)/2} \sum_{k=1}^N
      X_k(t_k-a) + B'(a, t).
   \end{equation}
   Here, $B'(a, t) = \int_{\Delta(a, t)} dW(s)$ and $W$
   is an $N$-parameter Brownian sheet in $\R$ independent of $B$,
   and all the processes on the right-hand side of
   (\ref{Eq:decomp}) are independent from one another.

   Thus,
   \begin{equation}\label{eq:BMLND}\begin{split}
      &\mathrm{Var} \left( \left. B(u) \,\right|\,
         B(t^1), \cdots, B(t^n) \right)\\
      & = \inf_{\alpha\in\R^n}
         \E\left[ \left( B(u) - \sum_{j=1}^n
         \alpha_j B(t^j)\right)^2
         \right]\\
      &\ge a^{N-1} \inf_{\alpha\in\R^{n}}
         \sum_{k=1}^N {\rm Var}
         \left( X_k(u_k-a) - \sum_{j=1}^n \alpha_j
         X_k(t_k^j-a)\right) \\
      &\ge a^{N-1} \sum_{k=1}^N {\rm Var}
         \left( \left. X_k(u_k-a) \, \right|\,
          X_k(t_k^1-a), \ldots, X_k(t_k^n-a) \right).
   \end{split}\end{equation}
   Therefore, (\ref{Eq:LND}) follows from (\ref{eq:BM1}).

   A simple modification of the preceding argument
   shows that (\ref{Eq:LND2}) follows from (\ref{eq:BM2});
   we omit the details.
\end{proof}

    We conclude this section with the following result.
\begin{lemma}\label{lem:linear}
   Let $n \ge 1$ be a fixed integer. Then for
   all distinct $t^1,\ldots,t^n\in (0, \infty)^N$, the random
   variables $B(t^1),\ldots,B(t^n)$ are linearly independent.
   \end{lemma}

\begin{proof}
   When all the coordinates of $t^1,\ldots,t^n$ are distinct, this
   follows from Proposition~\ref{Prop:LND}. In general, it suffices
   to show that for all constants $\alpha_1, \ldots, \alpha_n \in
   \R$, if
   $\mathrm{Var}(\sum_{j=1}^n \alpha_j B(t^j)) = 0$, then
   $\alpha_1 = \cdots = \alpha_n = 0$.
   The said variance is equal to
   $\int_{\R^N_+} (\sum_{j=1}^n \alpha_j \1_{[0,\, t^j]}(s))^2\, ds$
   ~\cite[p.~142]{Kh},
   which is assumed to be zero.
   Thus,
   $\sum_{j=1}^n \alpha_j \mathbf{1}_{[0,t^j]}(s)=0$ for a.e. $s \in
\R_+^N$,
   whence $\alpha_1=\cdots=\alpha_n=0$.
\end{proof}

\section{Analysis of Bridges} \label{Sec:Bridge}

For all $s\in\R^N_+$ we define the process $\{ B_{s}(t)
\}_{t\in\R^N_+}$ as
\begin{equation}\label{eq:bridge1}
   B_s(t) = B(t) - \prod_{j=1}^N \left(
   \frac{s_j\wedge t_j}{s_j}\right) B(s),
   \qquad{}^\forall t\in\R^N_+.
\end{equation}

In the case that $s$ has some coordinates that are zero we define
$0\div 0=1$ to ensure that the preceding is well-defined.
Clearly, $B_s(s)=0$ and $B_ {0}=B$. Thus, the process $B_s$ is a
realization of the sheet $B$ ``conditioned to be zero at time
$s$.'' Alternatively, $B_s(t)$ is the conditional least-squares
estimator of $B(t)$ given $B(s)$; i.e.,
\begin{equation}\label{eq:regression}
   B_s(t) = B(t) - \E \left[\left. B(t)\, \right|\, B(s) \right].
\end{equation}
Hence, for all fixed $s, t \in \R_+^N$, $B_s(t)$ is independent of
$B(s)$. It turns out that much more is true, viz.,

\begin{lemma}\label{lem:MP}
   Fix a partial order $\pi\in\Pi_N$ and a point
   $s\in\R^N_+$. Then,
   \begin{equation}
      \left\{ B_s(t)\right\}_{t\supi s}
      \text{ is independent of }\F_\pi(s).
   \end{equation}
\end{lemma}

\begin{proof}
   Because $B_s$ is a Gaussian process it suffices to check that
   if $t\supi  s\supi  u$, then $\E[B_s(t) B(u)]=0.$ Now
   \begin{equation}\begin{split}
      \E\left[ B(t) B(u) \right] & = \prod_{j\in\pi}u_j
      \cdot \prod_{j\not\in\pi} t_j,\\
      \E\left[ B(s) B(u) \right] & = \prod_{j\in\pi}u_j
      \cdot \prod_{j\not\in\pi} s_j.
   \end{split}\end{equation}
   The lemma follows
   because $\prod_{j=1}^N (s_j\wedge t_j)=\prod_{j\in\pi}s_j\cdot
   \prod_{j\not\in\pi}t_j$.
\end{proof}

Next we prove that the local dynamics of the bridge $\{
B_s(t)\}_{t\supi  s}$ are similar to those of the sheet $B$;
compare to Lemma~\ref{lem:Q(0)}.

\begin{lemma}\label{lem:varest}
   For each $a>0$, all
   partial orders $\pi\in\Pi_N$,
   and every $s$, $u$ and $v \in[a, b]^N$ that satisfy $s\lepi u,v$,
   \begin{equation}\label{Eq:compare}
      \frac{a^{N-1}}{2N} \|u-v\| \le
      \E\left[ \left( B_s(u)-B_s(v) \right)^2 \right] \le N
      b^{N-1}\|u-v\|.
   \end{equation}
\end{lemma}

\begin{proof}
   Thanks to (\ref{eq:regression}),
   \begin{equation}\label{eq:condexp}
      \E\left[ \left( B_s(u)-B_s(v) \right)^2 \right]
      = \mathrm{Var}\left( \left. B(u)-B(v)\,\right|\, B(s)\right).
   \end{equation}
   By (\ref{eq:projection}), this is bounded above by
   $\E[(B(u)-B(v))^2]$, which is at most
   $N\, b^{N-1}\|u-v\|$; cf.\@ Lemma~\ref{lem:Q(0)}.
   This proves the upper bound.

   For the lower bound, we apply (\ref{Eq:LND2}) in (\ref{eq:condexp})
   and obtain
   \begin{equation}\begin{split}
      \E\left[ \left( B_s(u)-B_s(v) \right)^2 \right]
         &\ge \frac{a^{N-1}}{2\sqrt{N}}\sum_{k=1}^N\min\left(
         |u_k-s_k|+|v_k-s_k| ~,~
         |u_k-v_k|\right)\\
      & = \frac{a^{N-1}}{2\sqrt{N}}\sum_{k=1}^N
         |u_k-v_k|,
   \end{split}\end{equation}
   owing to the triangle inequality. The lower bound follows.
\end{proof}

\begin{lemma}\label{lem:ballest}
   Fix two numbers $0<a<b<\infty$.
   Then there exists a finite constant $A_{\ref{eq:ballest}}>1$,
   which depends only on $(N,a,b)$, such that for all
   $s,t\in[a,b]^N$ and all $\e>0$,
   \begin{equation}\label{eq:ballest}
      A_{\ref{eq:ballest}}^{-1}
      \exp\left( - \frac{\e^2}{A_{\ref{eq:ballest}}\, \|s-t\|}
      \right)\cdot
      \frac{\e}{\|s-t\|^{1/2}} \le
      \P\left\{ |B_s(t)| \le\e \right\} \le A_{\ref{eq:ballest}}\,
      \frac{\e}{\|s-t\|^{1/2}}.
   \end{equation}
\end{lemma}

\begin{proof}
   We derive the upper bound first.
   \begin{equation}\begin{split}
      \P\left\{ |B_s(t)| \le \e\right\} & = \frac{1}{
      \sqrt{2\pi\mathrm{Var }B_s(t)}} \int_{-\e}^\e
      \exp\left( -\frac{z^2}{\mathrm{Var }B_s(t)} \right)\, dz\\
      & \le \e\sqrt{\frac{2}{\pi\mathrm{Var }B_s(t)}}.
   \end{split}\end{equation}
   For the $s$ and $t$ in question
   we can find a partial order $\pi \in \Pi_N$
   such that $s \lepi t$. Therefore, thanks to
   Lemmas~\ref{lem:Q(0)} and \ref{lem:varest},
   we have,
   \begin{equation}\label{eq:VarB}
      \frac{a^{N-1}}{2 N} \|s-t\| \le
      \mathrm{Var }B_s(t) \le b^{N-1}\|s-t\|.
   \end{equation}
   The upper bound follows. The lower bound is derived similarly.
\end{proof}

We end with a final elementary lemma on Gaussian ball-estimates:

\begin{lemma}\label{lem:gauss}
   Suppose $Y$ is a centered one-dimensional
   Gaussian random variable with variance $\s^2$.
   Let $\alpha$ and $\beta$ be two fixed positive numbers.
   Then, for all $x\in[-\alpha\s,+\alpha\s]$,
   \begin{equation}\displaystyle
      \P \{ |Y-x| \le \e\} \ge
      \begin{cases}
         e^{-\frac12 \alpha^2 -\alpha\beta} \P\{ |Y|\le\e\}, & \text{ if
            $\e\le \beta\s$},\\ \displaystyle
         \sqrt{\frac{2}{\pi}}\ \beta e^{-\s^2(\alpha+\beta)^2/2},
            & \text{ if $\e> \beta\s$}.
      \end{cases}
   \end{equation}
\end{lemma}

\begin{proof}
   Evidently,
   \begin{equation}\label{Eq:what}\begin{split}
      \P\{|Y-x|\le\e\} & = \int_{-\e}^{\e}
         e^{-(z+x)^2/(2\s^2)}\, \frac{dz}{\s\sqrt{2\pi}}\\
      & \ge e^{-\frac12 \alpha^2 -\alpha\e/\s} \int_{-\e}^\e
         e^{-z^2/(2\s^2)}\, \frac{dz}{\s\sqrt{2\pi}}.
   \end{split}\end{equation}
   When $\e\le \beta\s$, the result follows immediately;
   when $\e >\beta\s$, use $\int_{-\e}^\e \ge
   \int_{-\beta\s}^{\beta\s}$ in the first line of (\ref{Eq:what}), and then
   change variables $[w=z/\s]$ to  deduce the lemma.
\end{proof}

\section{Proof of Theorem~\ref{th:Kahane}}

\subsection{First Part}
We can first consider $F_n=F\cap[1/n,n]^N$, prove the theorem
with $F$ replaced by $F_n$, and then let $n\uparrow\infty$. This
shows that we might as well assume the following:
\begin{equation}\label{eq:Fab}
   F \subseteq [a,b]^N,\ \text{ where }0<a <b <\infty.
\end{equation}

For all $x\in\R^d$, $\e>0$, and $\mu\in\mathscr{P}(F)$, define
\begin{equation}\label{eq:L}
   l_\mu^\e(x) =\int \frac{\1_{\{ |B(s)-x|\le \e\}}}{
   (2\e)^d}\, \mu(ds).
\end{equation}
Also define $p_t$ to be the probability density function of
$B(t)$; i.e.,
\begin{equation}\label{eq:p}
   p_t(x) = \frac{e^{-\|x\|^2/(2 \prod_{j=1}^N t_j)}}{(2\pi
   \prod_{j=1}^N t_j)^{d/2}},\qquad{}^\forall x\in\R^d,\
   t\in\R^N_+.
\end{equation}
By Fatou's lemma,
\begin{equation}\begin{split}
   \E\left[ l_\mu^\e(x) \right]  & = (2\e)^{-d}
      \int_F \int_{\{|y-x|\le\e\}} p_s(y)\, dy\, \mu(ds)\\
   & \ge (1+o(1))\int p_s(x)\, \mu(ds)\quad(\e\to 0).
\end{split}\end{equation}
Thanks to (\ref{eq:Fab}), we can find a positive and finite
constant $A_{\ref{eq:EL}}=A_{\ref{eq:EL}}(N,d,a,b)$, such that
for all $\e\in(0,1)$, $\mu\in\mathscr{P}(F)$, and  $x\in\R^d$,
\begin{equation}\label{eq:EL}
   \E\left[ l_\mu^\e(x) \right]  \ge A_{\ref{eq:EL}}\exp\left(
   -\frac{\|x\|^2}{2 A_{\ref{eq:EL}}}\right).
\end{equation}

\begin{lemma}\label{lem:joint}
   Given (\ref{eq:Fab}) there exists a positive finite
   constant $A_{\ref{eq:joint}} =A_{\ref{eq:joint}} (d,N,a)$
   such that for all $x\in\R^d$,
   $s,t\in\R^N_+$, and $\e>0$,
   \begin{equation}\label{eq:joint}
      \P\left\{ |B(s)-x|\le \e ~,~ |B(t)-x|\le \e\right\}
      \le A_{\ref{eq:joint}}  \left( \frac{(2\e)^2}{\|t-s\|^{1/2}}
      \wedge 1\right)^d.
   \end{equation}
   Consequently,
   \begin{equation}\label{eq:joint1}
      \E\left[ \left( L^\e_\mu(x) \right)^2 \right]
      \le A_{\ref{eq:joint}} \ I_{d/2}(\mu).
   \end{equation}
\end{lemma}

\begin{proof}
   We will derive (\ref{eq:joint}); (\ref{eq:joint1}) follows
   from (\ref{eq:joint}) and a Fubini--Tonelli argument.
   We also note that because of
   the independence of the coordinates, it suffices to prove
(\ref{eq:joint}) when $d=1$.

   Define
   \begin{equation}\label{eq:c}
      C_{s,t} = \prod_{j=1}^N \left( \frac{s_j\wedge t_j}{s_j}\right).
   \end{equation}
   Note that $0\le C_{s,t}\le 1$.
   Then, recall (\ref{eq:bridge1}) and Lemma~\ref{lem:MP}
   to deduce that
   \begin{equation}\begin{split}
      &\P\left\{ |B(s)-x|\le \e ~,~ |B(t)-x|\le\e\right\}\\
      & = \P\left\{ |B(s)-x|\le \e ~,~ |B_s(t) + C_{s,t} B(s) -x|\le\e
\right\}\\
      & \le \P\left\{ |B(s)-x|\le \e\right\}\cdot
         \sup_{z\in\R} \P \left\{ \left| B_s(t)+z\right| \le \e \right\}.
   \end{split}\end{equation}
   Gaussian laws are unimodal, and this means that the supremum is
   achieved at $z=0$; i.e.,
   \begin{equation}
      \P\left\{ |B(s)-x|\le \e ~,~ |B(t)-x|\le\e\right\}
      \le \P\left\{ |B(s)-x|\le \e\right\}\cdot
      \P \left\{ \left| B_s(t)\right| \le \e \right\}.
   \end{equation}
   Thanks to (\ref{eq:p}) and (\ref{eq:Fab}),
   \begin{equation}
      \P\{ |B(s)-x|\le\e\}  = \int_{x-\e}^{x+\e}
       p_s(y)\, dy \le 2\e p_s(0) \le \e\sqrt{\frac{2}{\pi a^N}}.
   \end{equation}
   On the other hand, by Lemma~\ref{lem:ballest},
   $\P\{ |B_s(t)| \le \e \}  \le A_{\ref{eq:ballest}}
   \e \|s-t\|^{-1/2}$, whence the lemma.
\end{proof}

We are ready to derive half of Theorem~\ref{th:Kahane}.

\begin{proof}[Proof of Theorem~\ref{th:Kahane}: First Half]
   Thanks to (\ref{eq:EL}), Lemma~\ref{lem:joint}, and the
   Paley--Zygmund inequality [see, e.g., Kahane~\cite{Kahane85a}*{p.\@ 8}],
   for all $\mu\in\mathscr{P}(F)$
   and all $\e\in(0,1)$,
   \begin{equation}\label{eq:NN}\begin{split}
      \P\left\{ \mathrm{dist}\left(
         x , B(F) \right) <\e \right\} & \ge
         \P\left\{ l_\mu^\e(x) >0 \right\}  \ge
         \frac{\left( \E [l_\mu^\e(x)] \right)^2}{%
         \E\left[ \left( l_\mu^\e(x)\right)^2\right] }\\
      & \ge  \frac{A_{\ref{eq:joint}} \cdot A_{\ref{eq:EL}}}{%
         I_{d/2}(\mu)}\exp\left(
   -\frac{\|x\|^2}{A_{\ref{eq:EL}}}\right).
   \end{split}\end{equation}
   The constants on the right-hand side do not depend on
   $\e\in(0,1)$ or $\mu\in\mathscr{P}(F)$. Let
   $\e\to 0$ and optimize over $\mu\in\mathscr{P}(F)$
   to deduce from the path-continuity of $B$ that the probability
   of the event $\{x\in B(F)\}$
   is at least $A_{\ref{eq:joint}}  A_{\ref{eq:EL}}
   \exp(-\|x\|^2/A_{\ref{eq:EL}})
   \C_{d/2}(F)$. Integrate this bound to deduce that
   whenever $\C_{d/2}(F)>0$, the expected value of $\l_d(B(F))$ is
   positive.
\end{proof}

\begin{remark}
   As we mentioned in the Introduction, we can also use the Fourier
   analytic method of
   J.-P.~Kahane~\citelist{\ycite{Kahane85a}\ycite{Kahane85b}}
   to prove that $\C_{d/2}(F)>0
   $ implies $\l_d(B(F)) > 0$ a.s. The constructive proof in this
   paper makes it possible to control the value of
   $\E\{\l_d(B(F))\}$ in terms of the capacity $\C_{d/2}(F)$.
\end{remark}

\subsection{Second Part: Step 1}

We divide the proof into three steps. In this first step, we
derive the main technical inequality which is equation
(\ref{eq:nearly2}) below. Throughout this portion of the argument,
$\mu$ is an arbitrary probability measure on the fixed compact set
$F\subset\R^N_+$, and $x\in\R^d$ is some fixed spatial point. We
also choose and fix a partial order $\pi\in\Pi_N$ throughout.

Define $C_{s,t}$ by (\ref{eq:c}). Then,
\begin{equation}\begin{split}
   \E\left[\left. l_\mu^\e(x) \, \right|\, \Fpi (s)\right]
      & \ge \int_{t\supi s} \frac{\P\left\{ \left.
      |B(t)-x|\le\e \, \right|\, \Fpi(s) \right\}}{(2\e)^d}\,
      \mu(dt)\\
   & = \int_{t\supi s} \frac{\P\left\{ \left.
      |B_s (t) + C_{s,t}B(s) -x|\le\e \, \right|\,
      \Fpi(s) \right\}}{(2\e)^d}\, \mu(dt).
\end{split}\end{equation}
Now, as events, we have the obvious inclusion,
\begin{equation}\begin{split}
   &\left\{ |B(s)-x| < \frac{\e}{2}\right\}\cap
      \left\{ |B_s(t)-(1-C_{s,t}) x|\le \frac{\e}{2}\right\}\\
   &\ \subseteq \left\{ |B_s (t) + C_{s,t}B(s) -x|\le \e \right\}.
\end{split}\end{equation}
The preceding two displays together yield the following bound:
Almost surely on the event $\{|B(s)-x| < \e/2\}$,
\begin{equation}\label{eq:nearly}
   \E\left[\left. l_\mu^\e(x) \, \right|\, \Fpi (s)\right]
   \ge \int_{t\supi s} \frac{\P\left\{
   |B_s (t) - \left( 1- C_{s,t}
   \right) x |\le \frac{\e}{2}
   \right\}}{(2\e)^d}\, \mu(dt).
\end{equation}
[The conditioning can be removed thanks to Lemma~\ref{lem:MP}.]

Next, for all $s\lepi t$,  (\ref{eq:c}) implies
\begin{equation}\begin{split}
   1-C_{s,t} & = \prod_{j\not\in\pi} s_j^{-1} \cdot \left[
      \prod_{j\not\in\pi} s_j - \prod_{j\not\in\pi}
      t_j\right]\\
   & = \prod_{j\not\in\pi} s_j^{-1} \cdot
   \E\left[ \left( B_1(u)-B_1(v) \right) B_1(u) \right],
\end{split}\end{equation}
where $B_1$ denotes the first coordinate process of the Brownian
sheet $B$, and $u$ and $v$ are defined as follows: For all
$j\in\pi$, $u_j=v_j=1$ and for all $j\not\in\pi$, $u_j=s_j$ and
$v_j=t_j$. By Lemma~\ref{lem:Q(0)}, the
Cauchy--Bunyakovsky--Schwarz inequality, and (\ref{eq:Fab}),
\begin{equation}\begin{split}
   1-C_{s,t} & \le \prod_{j\not\in\pi} s_j^{-1} \cdot
      \sqrt{\E\left[ \left( B_1(u)-B_1(v) \right)^2\right] \cdot
      \E\left[ \left( B_1(u) \right)^2 \right]}\\
   & \le \prod_{j\not\in\pi} s_j^{-1} \cdot \sqrt{N} \, b^{(N-1)/2}
   \|u-v\|^{1/2} \prod_{j=1}^N u_j^{1/2}\\
   & \le \sqrt{N}\, b^{(N-1)/2} \|s-t\|^{1/2}.
\end{split}\end{equation}
%Here, $|\pi|$ denotes the cardinality of $\pi\in\Pi_N$.
In particular, we can find a positive and finite constant
$A_{\ref{eq:nearly1}}=A_{\ref{eq:nearly1}}(a,b,d,N)$ such that for
all $s,t\in[a,b]^N$, $1-C_{s,t}\le A_{\ref{eq:nearly1}}
\|s-t\|^{1/2}$. Plug this into (\ref{eq:nearly}) to obtain the
following: Almost surely on the event $\{|B(s)-x| < \e/2\}$,
\begin{equation}\label{eq:nearly1}\begin{split}
   &\E\left[\left. l_\mu^\e(x) \, \right|\, \Fpi (s)\right]\\
   &\ \ge (2\e)^{-d} \int_{t\supi s}\ \inf_{|z|\le
      A_{\ref{eq:nearly1}} \|s-t\|^{1/2}\cdot |x|}
      \P\left\{ |B_s (t)- z | < \frac{\e}{2} \right\}\, \mu(dt).
\end{split}\end{equation}
Hold $x\in\R^d$ fixed.  We can deduce from Lemmas~\ref{lem:gauss}
and~\ref{lem:ballest}, and equation (\ref{eq:VarB}), that there
exists a finite constant $A_{\ref{eq:nearly2}} =
A_{\ref{eq:nearly2}}(a,b,d,N,x)\in(0,1)$ such that whenever
$\e\in(0,A_{\ref{eq:nearly2}})$,
\begin{equation}\label{eq:nearly2}\begin{split}
   &\E\left[\left. l_\mu^\e(x) \, \right|\, \Fpi (s)\right]\\
   &\ \ge A_{\ref{eq:nearly2}} \int_{t\supi s}
      \left[ \frac{1}{\max\left( \e ~,~ \|t-s\|^{1/2}
      \right)} \right]^d \mu (dt)
      \cdot \1_{\{ |B(s)-x| <\e/2\}}.
\end{split}\end{equation}

\subsection{Step 2}
For the second portion of our proof, let us assume that $F$ has a
nonempty interior, and of course (\ref{eq:Fab}) is enforced as
well.

We will also make use of the fact that $F$ has a countable dense
subset. For simplicity, we assume it is a subset of the rational
numbers $\Q^N_+$. By continuity, the distance between $x$ and
$B(F)$ is less than $\e$ if and only if there exists a rational
time-point $t\in\Q^N_+\cap F$ such that
$\mathrm{dist}(x,B(t))<\e$. Moreover, the absolute-continuity of
the distribution of $B(s)$---for a given rational time-point
$s$---tells us that the latter happens with positive probability.
But it can happen also that with some positive probability
$\mathrm{dist}(x,B(t))\ge \e$.

In order to properly describe this last assertion, we let
$\partial\not\in\R^N_+$ denote a cemetery-point (in time), and
define $\Q^N_\partial=\Q^N_+\cup\{\partial\}$. Now enumerate all
rational time-points to deduce the existence of a
$(\Q^N_\partial\cap F)$-valued random variable $T_\e$ such that:
\begin{enumerate}
   \item $T_\e=\partial$ if and only if
         $\mathrm{dist}(x,B(F))\ge\e$;
   \item On the event $\{T_\e\neq\partial\}$,
         $T_\e \in F$ (a.s.), and
         $\mathrm{dist}(x,B(T_\e))<\e$.
\end{enumerate}
Because (\ref{eq:nearly2}) holds almost surely simultaneously for
all rational time-points $s$ and all partial orders $\pi\in\Pi_N$,
it follows that $\sup_{s\in\Q^N_+} | \E[  l_\mu^\e(x)\, |\,
\Fpi(s) ] |^2$ is bounded below by
\begin{equation}\label{eq:BB}
   \left( A_{\ref{eq:nearly2}}
   \int_{t\supi s} \left[ \frac{1}{\max\left( \e ~,~
   \|T_{\e/2} -t\|^{1/2}\right) } \right]^d\, \mu(dt)
   \right)^2 \cdot \1_{\{ T_{\e/2}\neq\partial\}}.
\end{equation}
So far, everything works for an arbitrary probability measure
$\mu$ on $F$. Now we describe a special choice for $\mu$. Namely,
we apply the preceding with $\mu$ replaced by
$\mu_\e\in\mathscr{P}(F)$, where
\begin{equation}
   \mu_\e (G) = \P\left\{ \left. T_{\e/2}  \in G
   \, \right|\, T_{\e/2}\neq\partial\right\},\qquad
   {}^\forall\text{ Borel sets }G\subseteq\R^N_+.
\end{equation}
Integrate (\ref{eq:BB}) [$d\P$] to conclude that
\begin{equation}\label{eq:WWW}\begin{split}
   & \E\left( \sup_{s\in\Q^N_+}
      \left| \E\left[ \left. l_{\mu_\e}^\e(x)\, \right|\,
      \Fpi(s) \right] \right|^2 \right)\\
   & \ge \E\left[ \left( A_{\ref{eq:nearly2}}
      \int_{t\supi s} \left[ \frac{1}{\max\left( \e ~,~
      \|T_{\e/2} -t\|^{1/2}\right) } \right]^d\, \mu_\e (dt)
      \right)^2 \cdot \1_{\{ T_{\e/2}\neq\partial\}}
      \right]\\
   & = A_{\ref{eq:nearly2}}^2 \int  \left(
      \int_{t\supi s} \left[ \frac{1}{\|s-t\|^{1/2}}\wedge
      \frac{1}{\e}
      \right]^d\, \mu_\e (dt)
      \right)^2\, \mu_\e (ds) \cdot \P\{T_{\e/2}\neq\partial\}\\
   & \ge A_{\ref{eq:nearly2}}^2 \left( \iint_{t\supi s}
      \left[  \frac{1}{\|s-t\|^{1/2}}\wedge \frac{1}{\e}
      \right]^d\, \mu_\e (dt) \, \mu_\e (ds) \right)^2
      \cdot \P\{T_{\e/2}\neq\partial\}\\
   & = A_{\ref{eq:nearly2}}^2 \mathscr{Q}_{\e,\pi}^2
      \cdot \P\{T_{\e/2}\neq\partial\}.
\end{split}\end{equation}
[In the fourth line, we have appealed to the
Cauchy--Bunyakovsky--Schwarz inequality.] On the other hand,
\begin{equation}\label{eq:WWWW}\begin{split}
    &\E\left( \sup_{s\in\Q^N_+}
      \left| \E\left[ \left. l_{\mu_\e}^\e(x)\, \right|\,
      \Fpi(s) \right] \right|^2 \right)\\
   &\ \le 4^N \sup_{s\in\Q^N_+} \E\left(
      \left| \E\left[ \left. l_{\mu_\e}^\e(x)\, \right|\,
      \Fpi(s) \right] \right|^2 \right)\\
   &\ = 4^N \E\left[ \left( l_{\mu_\e}^\e(x) \right)^2 \right]\\
   &\ \le A_{\ref{eq:WWWW}} \iint \left[  \frac{1}{\|s-t\|^{1/2}}\wedge
      \frac{1}{\e} \right]^d\, \mu_\e (dt) \, \mu_\e (ds)\\
   &\ = A_{\ref{eq:WWWW}} \mathscr{Q}_\e.
\end{split}\end{equation}
Brief justification: The first line follows from
Corollary~\ref{co:commutation}; and the third line follows from
(\ref{eq:joint}) and the Fubini--Tonelli theorem. We reemphasize
that the constants $A_{\ref{eq:nearly2}}$ and $A_{\ref{eq:WWWW}}$
do not depend on $\e$ or $\pi$. Add the preceding over all
$\pi\in\Pi_N$ to obtain
\begin{equation}\label{Eq:Square}
\begin{split}
   2^N A_{\ref{eq:WWWW}} \mathscr{Q}_\e &\ge A_{\ref{eq:nearly2}}^2
\sum_{\pi\in\Pi_N}
      \mathscr{Q}_{\e,\pi}^2 \cdot \P\{ T_{\e/2}\neq\partial\}\\
   & \ge 2^{1-N} A_{\ref{eq:nearly2}}^2 \left( \sum_{\pi\in\Pi_N}
\mathscr{Q}_{\e,\pi}
      \right)^2 \cdot \P\{ T_{\e/2}\neq\partial\}\\
   & \ge 2^{1-N} A_{\ref{eq:nearly2}}^2 \mathscr{Q}_\e^2 \cdot \P\{
T_{\e/2}\neq\partial\}.
\end{split}\end{equation}
Solve for the probability, using the fact that $\mathscr{Q}_\e$ is
strictly positive, to obtain
\begin{equation}
   \P\left\{ \mathrm{dist}(x,B(F))< \e/2\right\}
   \le \frac{2^{2N-1} A_{\ref{eq:WWWW}}}{A_{\ref{eq:nearly2}}^2
      \displaystyle\iint  \left[  \frac{1}{\|s-t\|^{1/2}}\wedge \frac 1\e
      \right]^d\, \mu_\e (dt) \, \mu_\e (ds)}
\end{equation}
Now, $\{\mu_\e\}_{\e\in(0,1)}$ is a collection of probability
measures on the compact set $F$; let $\mu_0$ denote any (weak)
limit-measure. Then, $\mu_0$ is also a probability measure on $F$,
and by the Fatou lemma and the path-continuity of $B$,
\begin{equation}\label{eq:good}
   \P\left\{ x\in B(F) \right\}
   \le \frac{2^{2N-1} A_{\ref{eq:WWWW}}}{A_{\ref{eq:nearly2}}^2
I_{d/2}(\mu_0)},
\end{equation}
and this is valid even if $\mu_0$ has infinite $\frac
d2$-dimensional energy as long as we interpret $1\div \infty$ as
zero.

\subsection{Step 3}
If $F$ has an interior, then (\ref{eq:good}) provides us with a
hitting estimate; we note once more that $A_{\ref{eq:nearly2}}$
and $A_{\ref{eq:WWWW}}$ of the latter equation depend only on
$(d,a,b,N,x)$. For a general compact set $F\subseteq[a,b]^N$, and
given $\eta\in(0,1)$, let $F^\eta$ denote the closed
$\eta$-enlargement of $F$. Equation~(\ref{eq:good}) provides us
with a positive finite constant $A_*=A_*(a,b,d,N,x)$ and a
probability measure $\mu^\eta$, on $F^\eta$, such that $\P\{ x\in
B(F^\eta)\} \le A_*/I_{d/2} (\mu^\eta)$. By the Fatou lemma and
weak compactness, we can find a probability measure $\mu_*$ on $F$
such that $\liminf_{\eta\to 0}I_{d/2}(\mu^\eta) \ge
I_{d/2}(\mu_*)$. Because $A_*$ does not depend on $\eta$,
path-continuity of $B$ shows that $\P\{ x\in B(F)\}\le
A_*/I_{d/2}(\mu_*) \le A_* \C_{d/2}(F)$. In particular, if $F$ has
zero $\frac d2$-dimensional capacity, then $\P\{ x\in B(F)\}=0$
for all $x$. A final appeal to the Fubini--Tonelli theorem
demonstrates that in this case, the expectation of $\l_d(B(F))$ is
zero. This completes our proof.

\section{Proof of Theorem \ref{th:level}}

Our proof of Theorem \ref{th:Kahane} contains the proof of Theorem
\ref{th:level}; cf.\@ (\ref{eq:NN}) and (\ref{eq:good}). However,
for the sake of future use we prove the following more general
result. It extends some results of
Kahane~\citelist{\ycite{Kahane72}\ycite{Kahane85b}} on stable
L\'evy processes.

\begin{proposition}\label{Prop:Kahane}
   Let $B = \{B(t)\}_{ t \in \R_+^N}$ denote the $(N,d)$
   Brownian sheet.
   Let $E \subset \R^d$ and $F \subset (0, \infty)^N$ be fixed
   Borel sets. Then the following are equivalent:
   \begin{enumerate}
       \item With positive probability, $F \cap B^{-1}(E) \ne
       \O$.
       \item With positive probability, $E \cap  B(F) \ne \O$.
       \item With positive probability, $\l_d ( E \ominus  B(F) ) >0$.
    \end{enumerate}
\end{proposition}

\begin{proof}
    Items (1) and (2) are manifestly equivalent. To prove $(2)
    \Leftrightarrow (3)$, we note that $(2)$ is equivalent to the
    following:
    \begin{equation}
       {}^\exists \delta > 0\ \text{ such that }\ \P \left\{ E \cap
       B\left(F \cap
       \left(\delta, \infty\right)^N\right) \ne \O \right\} >0.
    \end{equation}
    Hence, without loss of generality, we can assume that $ F \subset
    (\delta, \infty)^N$ for some $\delta > 0$. Fix a $\tau \in (0,
    \delta)^N$ so that $\tau \cle  t$ for all $ t \in F$. Define the
    random field $B^\tau = \{B^\tau (t)\}_{t \in \R_+^N}$ by
    \begin{equation}
    B^\tau(t) = B(\tau + t) - B(\tau), \quad  {}^\forall  t \in \R_+^N.
    \end{equation}
    Observe that
    \begin{equation}\label{Eq:48}
       E \cap  B(F) = \O \ \Longleftrightarrow\  B( \tau) \notin E
       \ominus  B^{\tau }(F- \tau).
    \end{equation}
    Because $B(\tau)$ is independent of the random Borel set $E \ominus
    B^\tau (F- \tau)$ and the distribution of  $B(\tau)$ is equivalent to
    $\l_d$, we have
    \begin{equation} \label{Eq:49}
       B( \tau) \notin E \ominus  B^{\tau}(F- \tau),\ \text{ a.s. }\
       \Longleftrightarrow\ \l_d \left( E \ominus
       B^{\tau}(F- \tau) \right) = 0, \ \text{ a.s. }
    \end{equation}
    Note that $ B^{\tau}(F- \tau) =  B(F)
    \ominus \{ B( \tau)\}$, so that the
    translation-invariance of the Lebesgue measure, (\ref{Eq:48}),
    and (\ref{Eq:49}) together imply that
    \begin{equation}
       E \cap  B(F) = \O,\ \text{ a.s.}  \ \Longleftrightarrow\
       \l_d \left( E \ominus  B(F)\right) = 0, \text{ a.s. }
    \end{equation}
    This proves the equivalence of (2) and (3), whence the proposition.
\end{proof}

\section{Proof of Theorem \ref{th:interior}}

Our proof of Theorem \ref{th:interior} relies on developing
moment-estimates for the local times of the Brownian sheet on $F$,
as well as a Fourier-analytic argument. Our argument is
closely-related to the methods
of~\ocite{Kaufman},~\ocite{Kahane85a},~\ocite{Mountford89},
and~\ocite{Xiao97}.

\subsection{First Reduction}
Without loss of generality, we can [and will] assume
that $F$ is compact. Otherwise, we can
consider a compact subset $F' \subseteq F$ such that
$\dimh  F' > \frac d2$; see~\ocite{Falconer}*{Theorem 4.10}.
Because $B(F')\subseteq B(F)$,
this proves that there is no harm in assuming that
(\ref{eq:Fab}) holds for some $0<a<b$. This
compactness assumption on $F$ is in force throughout this
section.

Because we have assumed that $\dimh  F > \frac d2$, we can choose
a $\gamma \in (0, 1)$ such that
\begin{equation}\label{eq:gamma}
   \dimh  F > \frac{\gamma  +  d} 2.
\end{equation}
Then by Frostman's lemma, there exists a probability measure $\mu$
on $F$ such that
\begin{equation}\label{Eq:Frostman}
   A_{\ref{Eq:Frostman}}
   = \sup_{s \in \R^N} \int \frac{\mu(dt)}{\left\|
   s - t \right\|^{(\gamma + d)/2}}<\infty.
\end{equation}
See~\ocite{Kahane85a}*{p.\@ 130} or~\ocite{Kh}*{p.\@ 517}.

\subsection{Second Reduction}
Fix some $c>0$ and define
\begin{equation}
   F_\ell = \left\{ t=(t_1,\ldots,t_N)\in F:\
   t_\ell = c\right\}\qquad\ell=1,\ldots,N.
\end{equation}
Suppose there exists $\ell\in\{1,\ldots,N\}$ such that
$\mu(F_\ell)>0$, where $\mu$ is the measure that satisfies
(\ref{Eq:Frostman}). By Frostman's lemma, the Hausdorff dimension
$F_\ell$ strictly greater than $(d/2)$. Identify $F_\ell$ with a
set in $\R^{N-1}_+$ (ignore the $\ell$th coordinate), and denote
the set in $\R^{N-1}_+$ by $F'_\ell$. The preceding development,
and Frostman's lemma, together prove that $\dimh(F_\ell')>(d/2)$.
It then suffices to prove that $\tilde{B}(F'_\ell)$ has an
interior point, where $\tilde{B}$ is $(N-1)$-parameter Brownian
sheet in $\R^d$. Therefore, we may assume---without loss of
generality---that the probability measure $\mu$ of
(\ref{Eq:Frostman}) has the following property: For all $c > 0$
and $\ell=1,\ldots,N$,
\begin{equation}\label{Eq:Fnot}
   \mu \left\{
   t= (t_1, \ldots, t_N)\in F:\ t_\ell = c \right\} = 0.
\end{equation}

Now consider the push-forward $\mu\circ B^{-1}$ of $\mu$ by
$B$. If $\mu\circ B^{-1} \ll \l_d$, then $B$ is said to have a
\emph{local time} on $F$. The local time $l_\mu(x)$ is defined as
the Radon--Nikod\'ym derivative $d\mu\circ B^{-1}/d \l_d$ at
$x\in\R^d$. Another way of writing this is this: If
$f:\R^d\to\R_+$ is Borel measurable, then with probability one,
\begin{equation}\label{eq:ODF}
   \int_F f(B(s))\, \mu(ds) = \int_{\R^d}
   f(x) l_\mu(x)\, dx.
\end{equation}

It is well known that (\ref{Eq:Frostman}) implies that $l_\mu$ is
in $L^2(\l_d)$ almost surely; cf.~\ocite{Geman}*{Theorem 22.1} or
\ocite{Kahane85a}*{Theorem 4, p.\@ 204}. In fact, $l_\mu$ is the
$L^2(\P\times\lambda_d)$-limit of $l^\e_\mu$ as $\e$ tends to $0$;
cf.\@ (\ref{eq:L}).

\subsection{Continuity in the Space-Variable}
Note that $B(F)$ is compact and $\{x: l_\mu(x)> 0\}$
is a subset of $B(F)$. Hence, in order to prove
that $B(F)$ has interior-points, it suffices to
demonstrate that $l_\mu(x)$ has
a version which is continuous in $x$~\citelist{\cite{Pitt78}*{p.\@
324}\cite{Geman}*{p.\@ 12}}. We are going to
have to do more to prove the uniform result for $\theta F$, but
for now, we concentrate on $\theta$ being equal to the
identity matrix.

\begin{theorem}\label{thm:LTcont}
   Let $F$ be a compact set in $\R^N_+$ that satisfies
   (\ref{eq:Fab}), and suppose $\mu\in\mathscr{P}(F)$ satisfies
   (\ref{Eq:Frostman}). Then for every even integer $n \ge 2$, there
   exists a finite constant $A_{\ref{Eq:m26}}$---that depends
   only on
   $(a,b,d,N, \gamma, A_{\ref{Eq:Frostman}}, n)$---such that
   \begin{equation}\label{Eq:m26}
         \E \left[ \sup_{\scriptstyle
         u,v\in\R^d:\atop\scriptstyle u\neq v}
         \frac{\left|
         l_\mu(u) - l_\mu(v)\right|^n}{|u-v|^{\lfloor
         n/(N+1)\rfloor \gamma}}
         \right] \le A_{\ref{Eq:m26}},
      \end{equation}
   where $\lfloor x \rfloor$ denotes 1 if $x \le 1$, and the largest
   integer $\le x$ if $x >1$. Consequently,
   there exists a version of $\{l_\mu(x)\}_{x\in\R^d}$ that is
   uniformly H\"older-continuous with index $\eta$
   for any $\eta$ that satisfies
   \begin{equation}
      0< \eta < \min\left(1 , \frac2 {N+1}
      \left(\dimh F- \frac d2 \right)\right),
   \end{equation}
\end{theorem}

Before proving Theorem~\ref{thm:LTcont}, we develop two technical
lemmas.

For the first lemma, define
\begin{equation}
   \widetilde{F}= \left\{\overline v = (v^1, \ldots, v^n) \in F^n; \,
   v^j_\ell= v^i_\ell \hbox{ for some } i \ne j \
   \hbox{ and }\ 1 \le \ell \le N\right\}.
\end{equation}

\begin{lemma}\label{lem:diagonal}
   If $\nu\in\mathscr{P}(F)$ satisfies (\ref{Eq:Fnot}), then the set
   $\widetilde{F}$ is $\nu^n$-null,
   where $\nu^n=\nu\times\cdots\times\nu$ ($n$ times).
\end{lemma}

\begin{proof}
   This follows from (\ref{Eq:Fnot}) and the Fubini--Tonelli theorem.
\end{proof}

\begin{lemma}\label{Lem:Cuzick}
   Let $\{Z_i\}_{i=1}^n$ be linearly-independent centered
   Gaussian variables. If $g:\R\to\R_+$ is Borel measurable, then
   \begin{equation}
      \int_{{\R}^n} g(v_1)\, e^{- \frac1 2 {\rm
      Var}\left(v\cdot Z\right)}\, dv
      = \frac{(2 \pi)^{(n-1)/2}} {Q^{1/2}}\
      \int_{-\infty}^\infty g(z/\sigma_1)\, e^{- \frac12 z^2}\, dz,
   \end{equation}
   where $\sigma_1^2 = {\rm Var}(Z_1\, |\, Z_2, \ldots, Z_n)$,
   and $Q = {\rm det\, Cov} (Z_1, \ldots, Z_n)$ denotes the determinant
   of the covariance matrix of $Z$.
\end{lemma}

\begin{proof}
   In the case that $g$ is bounded this follows
   from~\ocite{Cuzick82}*{Lemma 2}. To prove the general
   case, replace $g$ by $g\wedge k$ and let $k$ tend to infinity.
\end{proof}

We will use the following elementary formula to estimate the
determinant of the covariance matrix of a Gaussian vector $Z$:

\begin{equation}\label{Eq:Gausscovdet}
   {\rm det\, Cov} (Z_1, \ldots, Z_n) = {\rm Var}(Z_1) \prod_{j=2}^n
   {\rm Var} \left(Z_j \left| \left\{ Z_i\right\}_{i \le j-1}\right.
   \right).
\end{equation}

 We are ready to present the following.

\begin{proof}[Proof of Theorem~\ref{thm:LTcont}]
   By the Fourier inversion theorem,
   for every $x,y \in \R^d$, and all even integers $n \ge 2$,
   \begin{equation} \label{Eq:moment1}\begin{split}
      &\E \left[ \left( l_\mu(x) - l_\mu(x+y) \right)^n \right]\\
      & = (2 \pi)^{-nd} \int_{F^n}
         \int_{{\R}^{nd}} \prod_{j=1}^n
         \left[ e^{- i  u^j \cdot x} - e^{- i u^j\cdot (x+y) }\right]\\
      &\qquad\times
         e^{-\frac12 \mathrm{Var}\left( \sum_{j=1}^n u^j \cdot
         B(t^j)\right)}\, d\overline{u}\, \mu^n(d\overline{t}).
   \end{split}\end{equation}
   Here,
   $\overline{u} = ( u^1, \ldots, u^n)$, $\overline{t} = (t^1, \ldots,
   t^n)$, and for each $j$,
   $u^j$ and $t^j$ are respectively in $\R^d$ and
   $(0, \infty)^N$.
   The details that lead to
   (\ref{Eq:moment1}) are explained in~\ocite{Geman}*{Eq.\@ 25.7};
   see also~\ocite{Pitt78}.

Consider the non-decreasing function $\Lambda(u) = \inf\{1,
|u|^{\gamma}\}$ and the elementary inequality
\begin{equation}\label{Ieq:ele}
   | e^{i u} - 1 | \le 2 \Lambda(u), \quad \forall u \in\R.
\end{equation}
This is valid because $\gamma$ is in $(0,1)$. By the
triangle inequality,
\begin{equation}\label{Eq:tri}
\big| e^{- i  u^j \cdot x} - e^{- i u^j\cdot (x+y) }\big| \le
\sum_{k=1}^d \big|e^{-iu^j_k y^{\ }_k} - 1\big|.
\end{equation}
    By expanding the product in (\ref{Eq:moment1}), using
   (\ref{Eq:tri}) and (\ref{Ieq:ele}), we obtain
   \begin{equation} \label{Eq:moment2}\begin{split}
      &\E \left[ \left( l_\mu(x) - l_\mu(x+y)\right)^n \right]\\
      &\le (2 \pi)^{-nd}\,2^n\, {\sum}^\prime \int_{F^n}
         \underbrace{\int_{{\R}^{nd}} \prod_{j=1}^{n}
         \Lambda( u^j_{k_j}\, y^{\ }_{k_j}) \,
         e^{- \frac 1 2 {\rm Var} \left( \sum_{j=1}^n
          u^j \cdot B(t^j) \right)} \, d\overline{u}}\limits_{
         = \mathscr{M} =
         \mathscr{M}\left(\overline k, \overline t\right)}\,
         \mu^n(d\overline{t}).
   \end{split}\end{equation}
   Here, $\sum^\prime$ signifies the sum over all sequences
   $\overline k = (k_1,
   \cdots, k_n)  \in \{ 1, \ldots, d\}^n$.
   In accord with Lemma~\ref{lem:diagonal}, the outer integral in
   (\ref{Eq:moment2}) can be taken to be over $F^n\backslash
   \widetilde{F}$.

   Fix $\overline k = (k_1, \cdots, k_n)  \in \{ 1, \ldots, d\}^n$
   and $t^1, \ldots, t^n \in  F$, we proceed to
   estimate  the integral $\mathscr{M}$
    in (\ref{Eq:moment2}). We will assume
    that $t^j_\ell$
   $(1 \le j \le n, 1 \le \ell \le N)$ are
   distinct (Lemma~ \ref{lem:diagonal}).
   Lemma \ref{lem:linear} implies that the Gaussian random
   variables $\{ B_k(t^j);\, k = 1, \ldots, d,\ j = 1, \ldots, n\}$ are
   linearly independent.  Hence, by applying the generalized
   H\"older's inequality, Lemma~\ref{Lem:Cuzick}, and the independence
   of the coordinate-processes $B_1, \ldots, B_d$ of $B$,
   the quantity $\mathscr{M}$ is bounded above by
   \begin{equation} \label{Eq:m24}\begin{split}
      & \prod_{j=1}^{n}\left\{
         \int_{{\R}^{n d}} \Lambda^n( u^j_{k_j}\, y^{\ }_{k_j}) \,
         e^{- \frac1 2 {\rm Var} \left( \sum_{i=1}^n \sum_{\ell=1}^d
         u^i_\ell B_\ell(t^i)
         \right)}\, d\overline {u} \right\}^{1/n} \\
      & = \frac{(2\pi)^{(nd-1)/2}}{[{\rm
         det\, Cov}(B_1(t^1),\ldots, B_1(t^n) )]^{d/2}} \\
      &\quad\times
         \prod_{j=1}^n \left\{ \int_{-\infty}^\infty
         \Lambda^n\Big( \frac{y^{\ }_{k_j}} {\sigma_j
         ( \overline t)} \, z\Big)\,
         e^{-\frac12 z^2}\, dz \right\}^{1/n}.
   \end{split}\end{equation}
   Here,  $\s_j^2 (\overline t)$ is the conditional variance of
   $B_{k_j}(t^j)$ given $B_\ell(t^i) $ $ (\ell \ne k_j$ and $1\le i \le n$,
or $\ell =
   k_j$ and $i \ne j)$.

   Since $B_1, \ldots, B_d$ are i.i.d., we have
   \begin{equation}\label{eq:sigmaj}
      \s_j^2 \left( \overline t \right)
      = \mathrm{Var} \left( B_1(t^j)\, \left|\,
      \left\{ B_1(t^i) \right\}_{i\neq j}\right. \right).
   \end{equation}

For $j=n$, we use $\Lambda(u) \le 2 |u|^\gamma$ and Stirling's
formula to derive
\begin{equation}\label{Eq:sn}
\left\{\int_{-\infty}^\infty \Lambda^n\Big( \frac{y^{\ }_{k_n}}
{\sigma_n ( \overline t)} \, z\Big)\,
         e^{-\frac12 z^2}\, dz\right\}^{1/n} \le A_{\ref{Eq:sn}} n^\gamma\,
\Big(\frac{\|y\|} {\sigma_n ( \overline t)}\Big)^\gamma.
\end{equation}
Hence, it follows from (\ref{Eq:Gausscovdet}), (\ref{Eq:moment2}),
(\ref{Eq:m24}), (\ref{eq:sigmaj}) and (\ref{Eq:sn}) that
\begin{equation}\label{Eq:IntF}\begin{split}
   &   \max_{\overline k\in\{1,\ldots,d\}^n}\, \int_{F^n} \mathscr{M}
    \left( \overline k,\overline t
    \right) \mu^n(d \overline t) \le  \int_{F^n}
    \frac{A_{\ref{Eq:IntF}}\, \|y\|^\gamma  }
    {[{\rm det\, Cov}(B_1(t^1),\ldots, B_1(t^{n-1}) )]^{d/2}} \\
   &\qquad \quad \times\frac1  {\sigma_n^{d +\gamma}( \overline t)}\,
    \prod_{j=1}^{n-1} \left\{ \int_{-\infty}^\infty
    \Lambda^n\left( \frac{\|y\|} {\sigma_j (\overline t)}
    \, z\right)\, e^{-\frac12 z^2}\, dz
    \right\}^{1/n} \mu^n(d \overline t),
\end{split}\end{equation}
where $A_{\ref{Eq:IntF}}$ is a constant depending on $d, \gamma $
and $n$ only.  We can estimate the preceding integral iteratively by
integrating in the order $\mu(dt^n)$, $\mu(dt^{n-1}),$ $\ldots,
\mu(dt^1)$.

Let $t^1, \ldots, t^{n-1} \in F$ be fixed points such that
$t^j_\ell$ $(1 \le j \le n-1,\, 1 \le \ell \le N)$ are distinct.
We consider the integral
\begin{equation} \label{Eq:N}
\mathscr{N} = \int_F \frac1 {\sigma_n^{d+\gamma} \left( \overline
t \right)}\,  \prod_{j=1}^{n-1} \left\{ \int_{-\infty}^\infty
 \Lambda^n\Big( \frac{\|y\|} {\sigma_j (\overline t)} \, z\Big)\,
  e^{-\frac12 z^2}\, dz \right\}^{1/n} \mu(d t^n).
\end{equation}
It follows from Proposition \ref{Prop:LND} that for every $1\le j
   \le n$,
   \begin{equation}\label{Eq:ConVar2}
      \sigma_j^2 \left( \overline
      t \right) \ge \frac{a^{N-1}} 2 \sum_{\ell=1}^N \min_{i \ne j}
      \left|t^i_\ell - t^j_\ell\right|.
   \end{equation}

   In order to estimate the sum in (\ref{Eq:ConVar2}) as a function
   of $t^n$, we introduce $N$ permutations $\Gamma_1, \ldots,
   \Gamma_N$ of $\{1, \ldots, n-1\}$ such that for every $\ell = 1,
   \ldots, N$,
   \begin{equation}\label{Eq:order}
   t^{\Gamma_\ell(1)}_\ell < t^{\Gamma_\ell(2)}_\ell < \ldots <
   t^{\Gamma_\ell(n-1)}_\ell.
   \end{equation}
   For convenience, we denote $t^{\Gamma_\ell(0)}_\ell = a $ and
  $t^{\Gamma_\ell(n)}_\ell = b$ for all $1 \le \ell \le N$.

  For every $(i_1,\ldots, i_N) \in \{1, \ldots, n-1\}^N$, let
  $\tau_{i_1,\cdots, i_N} = (t^{\Gamma_1(i_1)}_1, \ldots,
  t^{\Gamma_N(i_N)}_N)$ be the ``center'' of the rectangle
  \begin{equation}
  I_{i_1,\cdots, i_N} = \prod_{\ell=1}^N
  \Big[t_\ell^{\Gamma_\ell(i_\ell)} - \frac1 2
  \big(t_\ell^{\Gamma_\ell(i_\ell)} -
  t_\ell^{\Gamma_\ell(i_\ell-1)}\big), \,
  t_\ell^{\Gamma_\ell(i_\ell)} + \frac 1 2
  \big(t_\ell^{\Gamma_\ell(i_\ell+1)} -
  t_\ell^{\Gamma_\ell(i_\ell)}\big)\Big)
  \end{equation}
with the convention that the left-end point of the interval is $a$
whenever $i_\ell=1$; and the interval is closed and its right-end
is $b$ whenever $i_\ell =n-1.$ Thus the rectangles
$\{I_{i_1,\cdots, i_N}\}$ form a partition of $[a, b]^N$.

For every $t^n \in F$, let $I_{i_1,\cdots, i_N}$ be the unique
rectangle containing $t^n$. Then (\ref{Eq:ConVar2}) yields the
following estimate:
   \begin{equation}\label{Eq:sigman}
     \sigma_n^2 \left( \overline t \right) \ge \frac{a^{N-1}} 2
      \sum_{\ell=1}^N \left|t^n_\ell -
      t_\ell^{\Gamma_\ell(i_\ell)}\right|
      \ge A_{\ref{Eq:sigman}} \|t^n - \tau_{i_1,\cdots, i_N}\|.
   \end{equation}

For every $j = 1, \ldots, n-1$, we say that $I_{i_1,\cdots, i_N}$
\emph{cannot see $t^j$ from direction $\ell$} if
   \begin{equation}
      t^j_\ell \notin \Big[ t_\ell^{\Gamma_\ell(i_\ell)} -
      \frac12 \big(
      t_\ell^{\Gamma_\ell(i_\ell)} - t_\ell^{\Gamma_\ell(i_\ell-1)}\big),
      \, t_\ell^{\Gamma_\ell(i_\ell)} +
      \frac12 \big(t_\ell^{\Gamma_\ell(i_\ell+1)}
      - t_\ell^{\Gamma_\ell(i_\ell)} \big) \Big].
   \end{equation}
We emphasize that if $I_{i_1,\cdots, i_N}$ cannot not see $t^j$
from all $N$ directions, then
\begin{equation}
 \big|t^j_\ell
- t^n_\ell \big| \ge \frac 1 2 \min_{i \ne j, n} \left|t^j_\ell -
t^i_\ell\right|\ \ \hbox{ for all }  1 \le \ell \le N.
\end{equation}
Thus $t^n$ does not contribute to the sum in (\ref{Eq:ConVar2}).
More precisely, the latter means that
\begin{equation}\label{Eq:sigmaj2}
 \sigma_j^2 \left( \overline t \right) \ge \frac{a^{N-1}} 4
\sum_{\ell=1}^N \min_{i \ne j, n} \left|t^j_\ell -
t^i_\ell\right|.
\end{equation}
The right hand side of (\ref{Eq:sigmaj2}) only depends on $t^1,
\ldots, t^{n-1}$, which will be denoted by $\tilde{\sigma}_j^2 (
\overline t )$. Hence we have
\begin{equation}\label{Eq:sigmaj3}
 \int_{-\infty}^\infty \Lambda^n\Big( \frac{\|y\|} {\sigma_j (\overline t)}
\,
z\Big)\,   e^{-\frac12 z^2}\, dz  \le \int_{-\infty}^\infty
\Lambda^n\Big( \frac{\|y\|} {\tilde{\sigma}_j (\overline t)} \,
z\Big)\, e^{-\frac12 z^2}\, dz.
\end{equation}

If $I_{i_1,\cdots, i_N}$ sees $t^j$ from a direction, then, except
in the special case $t^j = \tau_{i_1,\cdots, i_N}$,  it is
impossible to control $\sigma_j^2 \left( \overline t \right)$ from
below as  in (\ref{Eq:sigman}) and (\ref{Eq:sigmaj2}) [recall that
$B$ is not LND with respect to $\E (B(u) - B(v))^2$]. We say that
$t^j$ is a ``bad point'' for $I_{i_1,\cdots, i_N}$. In this case,
we use the inequality $\Lambda(u) \le 2$ to derive
\begin{equation}\label{Eq:sigmaj1}
 \int_{-\infty}^\infty \Lambda^n\Big( \frac{\|y\|} {\sigma_j (\overline t)}
\,
z\Big)\,  e^{-\frac12 z^2}\, dz  \le 2^n\sqrt{2 \pi}.
\end{equation}
It is important to note that, because of (\ref{Eq:order}), the
rectangle $I_{i_1,\cdots, i_N}$ can only have at most $N$ bad
points $t^j$ ($ 1 \le j \le n-1$), i.e., at most one in each
direction.

It follows from (\ref{Eq:sigman}), (\ref{Eq:sigmaj3}) and
(\ref{Eq:sigmaj1}) that
\begin{equation}\label{Eq:cubeI}
\begin{split}
   & \int_{I_{i_1,\cdots, i_N}} \frac1 {\sigma_n^{d+\gamma} \left(
    \overline t \right)}\,  \prod_{j=1}^{n-1} \left\{
    \int_{-\infty}^\infty
        \Lambda^n\Big( \frac{\|y\|} {\sigma_j (\overline t)} \, z\Big)\,
        e^{-\frac12 z^2}\, dz \right\}^{1/n} \, \mu(d t^n)\\
   &\le A\, \int_{I_{i_1,\cdots, i_N}} \frac {\mu(dt^n) } {\big\|t^n
    - \tau_{i_1,\cdots, i_N}\big\|^{(d +\gamma)/2}} \\
   &\quad \times \prod_{j
    \notin \Theta_{i_1\ldots i_N}}
    \left\{ \int_{-\infty}^\infty
        \Lambda^n\Big( \frac{\|y\|} {\tilde{\sigma}_j
    (\overline t)} \, z\Big)\,
    e^{-\frac12 z^2}\, dz \right\}^{1/n} \\
   & \le A_{\ref{Eq:cubeI}}\, \prod_{j \notin \Theta_{i_1\ldots
    i_N}} \left\{ \int_{-\infty}^\infty
    \Lambda^n\Big( \frac{\|y\|} {\tilde{\sigma}_j
    (\overline t)} \, z\Big)\,
    e^{-\frac12 z^2}\, dz \right\}^{1/n},
\end{split}
\end{equation}
where $\Theta_{i_1, \ldots, i_N} = \{1\le j \le n-1: \,  t^j\
\hbox{ is a bad point for }\ I_{i_1,\cdots, i_N}  \}$ and the last
inequality follows from (\ref{Eq:Frostman}). Recall that the
cardinality of $\Theta_{i_1, \ldots, i_N} \le N$ and
$\Theta_{i_1, \ldots, i_N}$ may be the same for different
$(i_1, \ldots, i_N)$.
%In the following, we find it more convenient to work with the
%complement of $J_{i_1, \ldots, i_N}$.

Summing (\ref{Eq:cubeI}) over all $(i_1,\ldots, i_N) \in \{ 1,
\ldots, n-1\}^N$ and regrouping $\Theta_{i_1\ldots i_N}$,
%the form $\{1, \ldots, n-1\}\backslash \Theta$,
we derive that the integral $\mathscr{N}$ is bounded above by
\begin{equation}\label{Eq:N2}
\begin{split}
& A_{\ref{Eq:cubeI}}\, \sum_{i_1,\ldots, i_N} \prod_{j \notin
\Theta_{i_1\ldots i_N}} \left\{ \int_{-\infty}^\infty
         \Lambda^n\Big( \frac{\|y\|} {\tilde{\sigma}_j (\overline t)} \,
z\Big)\,
         e^{-\frac12 z^2}\, dz \right\}^{1/n} \\
&\le A_{\ref{Eq:N2}}\, \sum_{\Theta} \,  \prod_{j\notin
\Theta} \left\{ \int_{-\infty}^\infty
         \Lambda^n\Big( \frac{\|y\|} {\tilde{\sigma}_j (\overline t)} \,
z\Big)\,
         e^{-\frac12 z^2}\, dz \right\}^{1/n},
\end{split}
\end{equation}
where the last summation is taken over all $\Theta \subset \{1,
\ldots, n-1\}$ with $\# \Theta \le N$ and $A_{\ref{Eq:N2}}$
depends on $(d, a, b, N,\gamma, n)$ only. Note that the number of
terms in the last sum is at most $(n-1)^N$.

Put (\ref{Eq:N2}) into (\ref{Eq:IntF}) to obtain
\begin{equation}\label{Eq:IntF2}\begin{split}
   &\max_{\overline k\in\{1,\ldots,d\}^n} \int_{F^n} \mathscr{M}
      \left( \overline k,\overline t  \right) \mu^n(d \overline t) \\
   &\le \sum_{\Theta} \, \int_{F^{n-1}} \frac {A_{\ref{Eq:IntF2}} \,
    \|y\|^{\gamma}} {[{\rm
        det\, Cov}(B_1(t^1),\ldots, B_1(t^{n-1}) )]^{d/2}}\\
   &\quad \times \prod_{j\notin \Theta} \left\{
    \int_{-\infty}^\infty
    \Lambda^n\Big( \frac{\| y\|} {\tilde{\sigma}_j
    (\overline t)} \, z\Big)\,
    e^{-\frac12 z^2}\, dz \right\}^{1/n}
    \mu^{n-1}(d \overline t).
\end{split}\end{equation}
We observe that:
\begin{enumerate}
    \item[(i)] Increasing the number of
       elements in $\Theta$ changes only the integrals  in
       (\ref{Eq:IntF2}) by a constant factor
       [recall (\ref{Eq:sigmaj1})]; and
    \item[(ii)] ${\rm det\, Cov}(B_1(t^1),\ldots, B_1(t^{n-1}))$ is
    symmetric in $t^1, \ldots, t^{n-1}$.
\end{enumerate}
Based on these observations we can deduce that
the following is valid uniformly for all
$\overline k \in \{1,\ldots,d\}^n$:
\begin{equation}\label{Eq:IntF3}\begin{split}
   \int_{F^n} \mathscr{M}
    \left( \overline k,\overline t  \right) \mu^n(d \overline t) &\le
    \int_{F^{n-1}} \frac {A_{\ref{Eq:IntF3}}\, \|y\|^{\gamma}} {[{\rm
    det\, Cov}(B_1(t^1),\ldots, B_1(t^{n-1}) )]^{d/2}}\\
   &\quad  \times \prod_{j=N+1}^{n-1} \left\{
    \int_{-\infty}^\infty
    \Lambda^n\Big( \frac{\| y\|} {\tilde{\sigma}_j
    (\overline t)} \, z\Big)\,
    e^{-\frac12 z^2}\, dz \right\}^{1/n}
    \mu^{n-1}(d \overline t).
\end{split}\end{equation}
Here, $A_{\ref{Eq:IntF3}}$ is a constant depending only on $(d, a,
b, N,\gamma, n)$, and the last product $\prod_{j=N+1}^{n-1}\cdots$
can be replaced by 1 if $n \le N+1$.

By repeating the preceding argument and integrating $\mu(dt^{n-1}),$
$\ldots, \mu(dt^1)$ iteratively and by (\ref{Eq:moment2}), we
obtain
\begin{equation} \label{Eq:moment3}
      \E \left[ \left( l_\mu(x) - l_\mu(x+y)\right)^n \right]
       \le  A_{\ref{Eq:moment3}}\,
       \|y\|^{\lfloor \frac{n} {N+1}\rfloor \gamma},
\end{equation}
where $A_{\ref{Eq:moment3}}$ is a constant depending on  $(d, a,
b, N,\gamma, n)$ only and $\|y\|^{\lfloor \frac{n} {N+1}\rfloor
\gamma}$ comes from the first $\lfloor \frac{n}
{N+1}\rfloor$-steps of integration.
   This, together with a multiparameter version of
   the Kolmogorov continuity theorem~\cite{Kh}*{Theorem 2.5.1, p.\@ 165}
   proves Equation (\ref{Eq:m26}), but where $u$ and $v$
   are restricted to a given compact set. It follows readily
   that we can construct a version of
   $\{l_\mu(x)\}_{x\in\R^d}$ that is
   H\"older-continuous with parameter $<\gamma$
   on all compact subsets of $\R^d$.
   But $x\mapsto l_\mu(x)$ is, by definition, a compact-support
   function because $\mu$ lives
   on the compact set $F\subset\R^N_+$, and $B$ is continuous.
   With a bit of measure theory, this completes our proof.
\end{proof}

\subsection{Continuity in the Rotation-Variable}

We hold $F$ and $\mu$ fixed as in the previous subsection. Let
$\mathscr{R}$ denote the collection of all $N$-by-$N$ rotation
matrices that leave $F$ in $\R^N_+$; i.e., $\theta\in\mathscr{R}$
if and only if $\theta$ is a rotation matrix such that $\theta
F\subset \R^N_+$. We endow all square matrices with their $\ell^2$
matrix-norm; i.e., for all $N$-by-$N$ matrices $M$,
\begin{equation}
   \| M \| = \sup_{x\in\R^N:\ \|x\|=1}
   \left( x\cdot Mx \right)^{1/2}
   = \sup_{x\in\R^N:\ \|x\|=1} \left\| Mx \right\|.
\end{equation}

Define,
\begin{equation}\label{eq:mu:theta}
   \mu_\theta (G) = \mu\left( \theta^{-1}G \right),
   \qquad{}^\forall \text{Borel sets }G\subset\R^N_+.
\end{equation}
Manifestly, $\mu_\theta$ is in $\mathscr{P}(F)$, and satisfies
(\ref{Eq:Frostman}) where $\mu$ is now replaced by $\mu_\theta$,
but the constant $A_{\ref{Eq:Frostman}}$ remains unchanged. Thanks
to Theorem~\ref{thm:LTcont}, $l_{\mu_\theta}$ is a.s.\@ H\"older
continuous for each $\theta\in\mathscr{R}$. We now prove that
there is a continuous version of $(x,\theta)\mapsto
l_{\mu_\theta}(x)$.

Throughout, we define
\begin{equation}\label{eq:eta}
   \tau(F) = \frac{\dimh F - \frac d2}{ 2\dimh F + N(d+1) +1}.
\end{equation}

\begin{theorem}\label{thm:LTcont2}
   Let $F$ be a compact set in $\R^N_+$ that satisfies
   (\ref{eq:Fab}), and suppose $\mu\in\mathscr{P}(F)$
   satisfies (\ref{Eq:Frostman}).
   For every even integer $n \ge 2$ and every
   $\delta<\tau(F)$,
   \begin{equation}\label{Eq:m27}
      \E \left[ \sup_{\scriptstyle \theta,\rho\in\mathscr{R}:\atop
      \scriptstyle \theta\neq\rho}\frac{
      \left| l_{\mu_\theta}(x) - l_{\mu_\rho}(x)\right|^n}{
      \|\theta-\rho\|^{n\delta}}
      \right] <\infty.
   \end{equation}
   Consequently, for all $0< \eta< \min(1 , \frac 2 {N+1}(\dimh F- \frac d2
))$
   and $0<\delta< \tau(F)$,
   there exists a version of $\{l_{\mu_\theta}
   (x)\}_{x\in\R^d; \theta\in\mathscr{R}}$ that is
   uniformly H\"older-continuous: In $x\in\R^d$ with index $\eta$,
   and in $\theta\in\mathscr{R}$ with index $\delta$.
\end{theorem}

\begin{proof}
   There exists a
   function $\psi:\R^d\to\R_+$ that has the following properties:
   \begin{itemize}
      \item There exists a finite constant $c_\psi$ such that
            for all $x,y\in\R^d$, $|\psi(x)-\psi(y)|\le c_\psi \|x-y\|$;
      \item $\psi(x)\ge 0$ and $\int_{\R^d}\psi(x)\, dx=1$;
      \item $\psi(x)>0$ if $\|x\|\le \frac12 $,
            whereas $\psi(x)=0$ if $\|x\|\ge 1$.
   \end{itemize}
   For all $\e>0$, define
   \begin{equation}
      \psi_\e (x) = \e^{-d} \psi(x/\e),\qquad{}^\forall x\in\R^d.
   \end{equation}

   First of all, note that for all $a\in\R^d$, $\theta\in\mathscr{R}$,
   and $\e>0$, the following holds a.s.:
   \begin{equation}\label{eq:EQ1}\begin{split}
      \left|\int_F \psi_\e \left( B(s)-a \right)\, {\mu_\theta}(ds)
         - l_{\mu_\theta}(a)
         \right| & = \left| \int_{\R^d} \psi_\e (x-a)\,
         l_{\mu_\theta}(x)\, dx  - l_{\mu_\theta}(a) \right|\\
      & \le \int_{\R^d} \psi_\e (x-a) \left| l_{\mu_\theta}(x)
         - l_{\mu_\theta}(a) \right|\, dx\\
      & \le \sup_{x\in\R^d: \,
         \| x-a\|\le \e} \left| l_{\mu_\theta}(x)
            -l_{\mu_\theta}(a)\right|\\
      & = \Omega_{\theta}(\e).
   \end{split}\end{equation}
   [Justification: The first line follows from (\ref{eq:ODF});
   second from the fact that $\psi_\e$ integrates to one;
   and third from the fact that $\psi_\e$ is supported on
   the centered ball of radius $\e$.]
   Furthermore, $\psi_\e$ is Lipschitz-continuous with Lipschitz-constant
   $c_\psi \e^{-(d+1)}$. Therefore, for all $\theta,\rho\in\mathscr{R}$,
   all $a\in\R^d$, and all $\e>0$, with probability one,
   \begin{equation}\label{eq:EQ2}\begin{split}
      & \left| \int_F \psi_\e \left( B(s)-a\right)\, {\mu_\theta}(ds)
         - \int_F \psi_\e \left( B(s)-a\right)\, \mu_\rho (ds)\right|\\
      & \ = \left| \int_F \left[
         \psi_\e \left( B(\theta s)-a\right) -
         \psi_\e \left( B(\rho s)-a\right)
         \right] \, \mu(ds) \right|\\
      & \ \le c_\psi \e^{-(d+1)} \sup_{s\in F} \left|
         B(\theta s) - B(\rho s) \right|.
   \end{split}\end{equation}
   Combine (\ref{eq:EQ1}) and (\ref{eq:EQ2}) to deduce that a.s.:
   \begin{equation}\label{eq:EQ3}
      \left| l_{\mu_\theta}(a) - l_{\mu_\rho}(a)
      \right| \le  \Omega_\theta(\e) + \Omega_\rho(\e)
      + c_\psi \e^{-(d+1)} \sup_{s\in F} \left|
      B(\theta s) - B(\rho s) \right|.
   \end{equation}
   One can use this to directly construct a
   continuous version of these
   local times. However, we will outline a more standard approach.

   By continuity (Theorem~\ref{thm:LTcont}), (\ref{eq:EQ3}) holds
   simultaneously for all $a\in\R^d$. Therefore, by
   Minkowski's inequality, for all even integers
   $n \ge 2$,
   \begin{equation}\begin{split}
      &\left( \E\left[ \sup_{a\in\R^d} \left|
         l_{\mu_\theta}(a) - l_{\mu_\rho}(a)
         \right|^n \right] \right)^{1/n}\\
      & \le  \left\|
         \Omega_\theta(\e)\right\|_{L^n(\P)} +
         \left\| \Omega_\rho(\e)\right\|_{L^n(\P)}\\
      &\qquad + c_\psi \e^{-(d+1)} \left( \E\left[ \sup_{s\in F} \left|
         B(\theta s) - B(\rho s) \right|^n \right] \right)^{1/n}.
   \end{split}\end{equation}
   Fix some positive $\eta<\min(1, \frac 2 {N+1}(\dimh F-\frac d2))$ to see
that
   the first two terms are each bounded above by
   $A_{\ref{Eq:m26}}\e^\eta $; see (\ref{Eq:m26}).
   Because $F$ is compact, standard Kolmogorov-continuity
   estimates show that
   the third term is at most a universal constant [depending only
   on $(d,N,a,b)$] times
   $\e^{-(d+1)}\|\theta-\rho\|^{1/2}$; for example,
   see Exercise 7 of~\ocite{Kh}*{p.~176}. Optimize the resulting inequality
   over all $\e$ to obtain (\ref{Eq:m27}). The remainder of the proof
follows
   from a multiparameter version of the Kolmogorov continuity theorem.
\end{proof}

\subsection{The Remainder of the Proof of Theorem~\ref{th:interior}}

We are ready to assemble the pieces that complete the proof of
Theorem~\ref{th:interior}. Throughout, we may, and will, assume
that $\{l_{\mu_\theta}(x)\}_{x\in\R^d; \theta\in\mathscr{R}}$ is
continuous (Theorem~\ref{thm:LTcont2}).

According to (\ref{eq:ODF}), we have
\begin{equation}
   \int_{\R^d} l_{\mu_\theta}(x)\, dx = 1,
   \qquad{}^\forall\theta\in\mathscr{R}.
\end{equation}
This uses the continuity of local times and $B$, as well as the
compactness of $F$. Note that we have stopped writing ``a.s.''
because from now on, there is  only one null-set left, and so it
can be ignored.

Continuity insures that for every $\theta\in\mathscr{R}$ there
exists an open ball $J_\theta\subset\R^d$ such that for all $x\in
J_\theta$, $l_{\mu_\theta}(x)>0$. This is enough to prove that for
all $\theta\in\mathscr{R}$, $B(\theta F)$ has interior-points: Any
$x\in J_\theta$ is an interior-point of $B(\theta F)$.

In order to prove the stronger assertion of the theorem, we need
to refine the $J_\theta$'s slightly.

Due to continuity, for every $\theta\in\mathscr{R}$ we can find an
open ball $K_\theta\subseteq J_\theta \subset \R^d$ and an open
ball $V_\theta\subset\mathscr{R}$ such that for all $\rho\in
V_\theta$ and all $x\in K_\theta$, $l_{\mu_\rho}(x)>0$. Now
$\{V_\theta\}_{\theta\in\mathscr{R}}$ is an open cover of
$\mathscr{R}$, where the latter is viewed as a closed subset of
$\mathfrak{H}_N$---the rotation group acting on $\R^N$. Because
$\mathfrak{H}_N$ is compact~\cite{P}*{Section 65, p.\@ 489}, so is
$\mathscr{R}$. It follows that there is a finite subcover
$\{V_{\theta(j)}\}_{j=1}^m$ of $\mathscr{R}$; it has the property
that for every $\rho\in V_{\theta(j)}$ and all $x\in
K_{\theta(j)}$, $l_{\mu_\rho}(x)>0$. Let $\zeta_j$ denote the
midpoint of the interval $K_{\theta(j)}$ to deduce the theorem.

\section{An Arithmetic Property of Brownian Motion}

We conclude this paper by proving an arithmetic result about
Brownian motion. Henceforth, $\{X(t)\}_{t\ge 0}$ denotes
$d$-dimensional Brownian motion, and $F$ a fixed compact subset of
$\R_+$.

Choose and fix an integer $N\ge 1$, and $N$ nonzero real numbers
$r_1,\ldots,r_N$, once and for all. Define the [inhomogeneous]
\emph{$N$-fold Brownian self-intersection field}~\cite{Wolpert}:
\begin{equation}
   S (t) = \sum_{j=1}^N r_j X(t_j),\qquad
   {}^\forall t\in\R^N_+.
\end{equation}

Next is a refinement to Theorem 1 of~\ocite{Mountford88}; see also
~\ocite{Kaufman79}.

\begin{theorem}\label{th:Mountford}
   Suppose $G_1,\ldots,G_N$ are compact subsets of
   $\R_+$, and let $G=G_1\times\cdots \times G_N$.
   If $\dimh G>\frac d2$, then $S(G)$
   a.s.~has interior-points. Moreover,
   $\C_{d/2}(G)=0$ if and only if $S(G)$ a.s.~has zero Lebesgue measure.
\end{theorem}

Before we prove this, we make some observations.
\begin{remark}
   Note the elementary bounds,
   \begin{equation}
      \sum_{j=1}^N \dimh (G_j)
      \le \dimh (G_1\times\cdots\times G_N) \le
      \sum_{j=1}^N \dimp (G_j).
   \end{equation}
   where $\dimp$ denotes the packing dimension.
   Therefore, Theorem~\ref{th:Mountford} implies
   that:
   \begin{itemize}
   \item If $\dimh (F) > \frac{d}{2N}$, then
      $X(F)\oplus\cdots\oplus X(F)$ contains interior-points a.s.
   \item If $\dimp (F) < \frac{d}{2N}$, then
      $X(F)\oplus\cdots\oplus X(F)$ is Lebesgue-null a.s.
   \end{itemize}
   The first item is a minor generalization
   of Theorems 2--4 of~\ocite{Mountford88};
   the second item slightly improves upon Comment (2)
   of~\ocite{Mountford88}*{p.\@ 459} who derives this assertion
   with upper Minkowski dimension in place of packing dimension.
\end{remark}

To prove Theorem~\ref{th:Mountford}
we may---and will---assume without loss of
generality that
$G_i \subset [a_i, b_i]$ and
\begin{equation}
   0 < a_1 < b_1 < a_2 < b_2< \cdots < a_N < b_N.
\end{equation}
Similar reductions have been made earlier
by~\ocite{Mountford88} and~\ocite{Kaufman79}. To
simplify the formulation of Proposition~\ref{prop:SI}, we
assume further that for all relevant integers $i$ and $j$,
$(b_i- a_i) \le (a_{j+1}- b_j)$.

Thanks to Lemma~\ref{lem:LND}, we can deduce

\begin{proposition}\label{prop:SI}
Suppose $\{G_i\}_{i=1}^N$ are compact subsets of $\R_+$
that satisfy the preceding conditions. Then
the process $S$ is sectorially LND on
$G=G_1\times\cdots \times G_N$. In fact, for all
$u,t^1,\ldots,t^n\in\R^N_+$,
   \begin{equation}\label{Eq:Arith}
      \mathrm{Var}\left( S_1(u) \, \left|\,
      S_1(t^1) , \ldots , S_1(t^n) \right. \right)
      \ge \frac{\min_{1\le k\le N}r_k^2}{2N}\sum_{k=1}^N
      \min_{1\le j\le n} \left| u_k - t^j_k \right|.
   \end{equation}
\end{proposition}

\begin{proof}\ Let $u = (u_1, \ldots, u_N)$. It follows from
Lemma~\ref{lem:LND} that for each $1 \le k \le N$
\begin{equation}
\begin{split}
\mathrm{Var}\left( S_1(u) \, \left|\,
      S_1(t^1) , \ldots , S_1(t^n) \right. \right)&\ge \mathrm{Var}\left(
      X(u_k) \left|\, X(t^j_\ell),\, \forall j,\, \ell;
      X(u_i), i \ne k \right. \right)\\
      &\ge \frac {r_k^2} 2\, \min_{1\le j\le n} \left| u_k - t^j_k \right|.
\end{split}
\end{equation}
Summing over $k$ yields (\ref{Eq:Arith}).
\end{proof}

\begin{proof}[Proof of Theorem~\ref{th:Mountford}]
   We go through the proof of Theorem~\ref{thm:LTcont},
   but use Proposition~\ref{prop:SI} in place of
   Proposition~\ref{Prop:LND} everywhere. This readily proves
   that when $\dimh G>\frac{d}{2}$ $S$ has a continuous local time
   on $G$. Therefore, $S(G)$ has  interior-points almost surely.

   For the capacity condition, we simply compare $S$
   to an additive Brownian motion. This is achieved
   by combining the proof of Theorem
   6.1 of~\ocite{KX03}
   with Theorem 4.2 in the same paper.
\end{proof}

\section{A Final Remark}
Consider an arbitrary centered Gaussian random field
$\{G(t)\}_{t\in\R^N_+}$,  a compact set $F$ in $\R^N_+$ that
satisfies (\ref{eq:Fab}) and $\mu\in\mathscr{P}(F)$ satisfies
(\ref{Eq:Frostman}) with $\gamma\in(0,1)$. An inspection of our
proof of Theorem~\ref{thm:LTcont} shows that everything up to and
including (\ref{Eq:N}) is valid as long as $G$ has i.i.d.\@
coordinate-processes. The rest of proof depends crucially on
whether one can derive an appropriate upper bound for
(\ref{Eq:N}).

If we apply the inequality $\Lambda(u) \le 2 |u|^\gamma$ to
(\ref{Eq:m24}) for all $j = 1, \ldots, n$, then for all
$x,y\in\R^d$ and all even integers $n\ge 2$ we can find a finite
constant $A_{\ref{eq:G:m24}} = A_{\ref{eq:G:m24}}(d,a,b,\gamma,
n)$ such that
\begin{equation}\label{eq:G:m24}\begin{split}
   & \E\left[ \left( l_\mu(x)-l_\mu(x+y) \right)^n \right]\\
   & \le A_{\ref{eq:G:m24}}\, \|y\|^{n\gamma}
      \int_{F^n}\frac{\prod_{j=1}^n\left[ \mathrm{Var} \left(
      G_1 (t^j) \, \left|\, \left\{G_1 (t^i)\right\}_{i\neq j}
      \right. \right) \right]^{-\gamma}}{ \left[
      \mathrm{det\, Cov}\left( G_1 (t^1),\ldots, G_1 (t^n) \right)
      \right]^{d/2} }  \, \mu^n(d\overline t).
\end{split}\end{equation}

Hence we have derived the following convenient result.

\begin{theorem}
   Let $\{G(t)\}_{t\in\R^N_+}$ be an $N$-parameter
   centered Gaussian process in $\R^d$ that
   has i.i.d.\@ continuous
   coordinate-processes. Choose and fix
   a compact set $F\subset \R_+^N$ such that
   $\sup_{t\in F} \mathrm{Var}(G_1(t)) <\infty.$
   Assume that there exists $\mu\in\mathscr{P}(F)$, an even
   integer $n > N$, and some $\gamma\in(\frac Nn,1)$ such that
   \begin{equation}
      \int_{F^n} \frac{\prod_{j=1}^n\left[ \mathrm{Var} \left(
      G_1 (t^j) \, \left|\, \left\{G_1 (t^i)\right\}_{i\neq j}
      \right. \right) \right]^{-\gamma}}{ \left[
      \mathrm{det\, Cov}\left( G_1 (t^1),\ldots, G_1 (t^n) \right)
      \right]^{d/2} }  \, \mu^n(d\overline t) <\infty.
   \end{equation}
   Then $\{l_\mu(x)\}_{x\in\R^d}$ has a modification
   which is H\"older continuous
   of any order $<\gamma$. Consequently,
   $G(F)$ has interior-points almost surely.
\end{theorem}

\end{document}